\documentclass[	]{elsarticle}
\pdfoutput=1

\newcommand{\mb}[1]{\mathbf{#1}}
\usepackage{amsmath,amsthm,amssymb,mathrsfs}

\newcommand{\norm}[1]{\| #1 \|}

\newcommand{\ud}{\mathrm{d}}
\newcommand{\oh}{\frac{1}{2}}
\usepackage{graphicx}
\usepackage{epstopdf}
\usepackage{subfig}
\usepackage{stackengine} 
\stackMath
\usepackage[usenames,dvipsnames]{color}
\usepackage{textcomp}
\newcommand{\edone}[1]{\textcolor{black}{#1}}
\newcommand{\edtwo}[1]{\textcolor{black}{#1}}
\newcommand{\edthree}[1]{\textcolor{black}{#1}}

\graphicspath{{figures/}{.}}
\usepackage{hyperref}

\journal{Journal of Computational Physics}

\begin{document}

\begin{frontmatter}

\title{Lubricated Immersed Boundary Method in Two Dimensions}


\author[mymainaddress]{Thomas G. Fai\corref{mycorrespondingauthor}}

\cortext[mycorrespondingauthor]{Corresponding author}
\ead{tfai@seas.harvard.edu}

\address[mymainaddress]{Paulson School of Engineering and Applied Sciences, Harvard University\\ 29 Oxford St., Cambridge, MA 02138}

\author[mymainaddress,berkeley]{Chris H. Rycroft}
\ead{chr@seas.harvard.edu}

\address[berkeley]{Mathematics Group, Lawrence Berkeley Laboratory\\ 1 Cyclotron Rd., Berkeley, CA 94720}

\begin{abstract}
Many biological examples of fluid-structure interaction, including the transit of red blood cells through the narrow slits in the spleen and the intracellular trafficking of vesicles into dendritic spines, involve the near-contact of elastic structures separated by thin layers of fluid. Motivated by such problems, we introduce an immersed boundary method that uses elements of lubrication theory to resolve thin fluid layers between immersed boundaries. We demonstrate 2\textsuperscript{nd}-order accurate convergence for simple two-dimensional flows with known exact solutions to showcase the increased accuracy of this method compared to the standard immersed boundary method. Motivated by the phenomenon of wall-induced migration, we apply the lubricated immersed boundary method to simulate an elastic vesicle near a wall in shear flow. We also simulate the dynamics of a vesicle traveling through a narrow channel and observe the ability of the lubricated method to capture the vesicle motion on relatively coarse fluid grids.
\end{abstract}


\end{frontmatter}


\section{Introduction}
The immersed boundary method is a widely-used numerical method for fluid-structure interaction that has been applied to problems including blood clotting \cite{crowl}, osmotic swelling due to thermal fluctuations \cite{wu2015}, sperm locomotion through viscoelastic fluids \cite{sperm}, and insect flight \cite{millercaf}. Several extensions of the immersed boundary method have been developed to incorporate realistic structural properties such as added mass \cite{kim2007penalty}, permeability \cite{Kim2006,Kim2012}, intrinsic twist \cite{slim}, nonlinear constitutive laws \cite{griffith2012hybrid}, and variable viscosity and density \cite{vc_ib,vc_ib2}. The chief idea of the immersed boundary method is to use an Eulerian description of the fluid and a Lagrangian description of the structure, while coupling these descriptions through integral operators involving delta functions. In the continuum formulation, these are the usual Dirac delta functions, and in practice much of the efficiency and accuracy of the immersed boundary method depends on how these Dirac delta functions are regularized \cite{acta,moridf}.

One of the features of the immersed boundary method is that it is formulated in a manner that avoids the contact problem; since all structures move in the same global velocity field, they cannot in principle cross themselves or one another \cite{lim2010dynamics}. Though it is useful for many purposes, advecting all structures in the same regularized velocity field presents challenges when multiple immersed boundaries are in near contact, which is typical in biological phenomena such as red blood cell motion through the microcirculation \cite{Pivkin2016} and membrane receptor trafficking in dendrites \cite{Kusters2014}. In simulations using fixed fluid grids, the flow in the thin fluid layer between immersed boundaries may not be sufficiently resolved, resulting in boundaries that are effectively stuck together. This is sometimes dealt with by the addition of repulsive forces, as done by Krishnan et al.~to  prevent ``unphysical overlap'' of rigid particles \cite{krishnan2017fully}, and by Lim et al.~in their simulations of flagellar bundling in \emph{E. Coli} to ``prevent filaments from crossing each other'' \cite{lim2012}. We wish to overcome this difficulty by directly using the governing equations instead.

The immersed boundary method takes the viewpoint that the entire domain, including any immersed elastic structures, is filled with fluid having a continuously varying velocity field, without large jumps in the velocities of structures that are nearby relative to the grid spacing $\Delta x$. 

Given the difficulties observed above, we propose a modification to the immersed boundary method that makes use of lubrication theory to resolve the flow in thin layers. Unlike the standard immersed boundary method, which becomes less accurate as two boundaries approach each other, the lubrication approximation actually \emph{improves} as the gap size decreases. The grid spacing $\Delta x$ provides a natural cutoff at which to apply the lubrication approximation. In the proposed method, the lubrication approximation is automatically turned off when the gap size becomes larger than the grid spacing, as described in detail later on.

We illustrate this method by applying it to two test problems of increasing complexity for which the lubricated immersed boundary method has significantly smaller errors than the standard method. The primary additional cost of the lubricated immersed boundary method involves the construction of a height function that gives the distance between the immersed boundaries. We describe a method to compute this height function accurately and efficiently. Finally, we apply the lubricated immersed boundary method to study the motion of elastic vesicles near walls and observe the significantly increased accuracy of simulations of wall-induced migration and channel flow on coarse fluid grids.

Our work builds on several previous efforts to merge direct numerical simulation with lubrication theory. In the context of particle suspensions, Nguyen and Ladd and Janoschek et al.~worked out lubrication corrections to lattice-Boltzmann simulations of particle suspensions \cite{nguyen2002lubrication,janoschek2013accurate}, while Seto et al.~performed simulations of hard spheres with regularized lubrication forces \cite{seto2013discontinuous}. Thomas et al.~performed a multiscale computation of a fluid drop interacting with a wall in which lubrication theory was used to replace the no-slip boundary condition by a specified wall-shear stress \cite{thomas2010multiscale}. Our approach is most closely related to this last work, since we also combine a grid-based fluid solver with analytical results from lubrication theory to model subgrid effects. Moreover, there are common features in terms of the mathematical formulation; we also use a height function to describe the size of the lubricating layer between boundaries, and relate the flow in the gap to the stress on the boundaries. However, whereas Thomas et al.~solve an evolution equation for the height and flux through the gap, the solution to which gives a boundary condition for the fluid, our boundaries are defined in a Lagrangian manner with the height function constructed based on an explicit representation of the boundaries as piecewise polynomial curves. In our approach, we use lubrication theory to modify the advection velocity. Further, whereas Thomas et al.~considered a multiphase fluid in which a drop of fluid having no surface elasticity and with ten times greater density and viscosity than the ambient fluid falls down a slope, in the present study we restrict attention to elastic structures immersed in fluids with constant density and viscosity.


We note that an alternative approach for ensuring that the fluid dynamics is sufficiently resolved throughout the domain is to use adaptive mesh refinement (AMR). The idea of AMR is to use a hierarchy of fluid grids with different levels of refinement, with the finer grids placed in specific regions of space to obtain a desired accuracy \cite{berger1984adaptive}. This allows for highly-resolved computations at drastically reduced computational cost compared to uniform grids. However, implementing AMR requires additional data structures to keep track of the grid hierarchy, and typical implementations of AMR require the user to specify a finest grid level so that issues with near-contact are not definitively resolved. This finest grid size can also lead to severe timestep restrictions for numerical stability. The advantage of using a subgrid model, as we do here through lubrication theory, is that the asymptotic results hold for infinitesimally small gaps, and does not introduce additional timestep constraints. Further, relatively little additional machinery is required beyond what is already needed for uniform grid simulations; as we will discuss later on, the main extra cost of the lubricated immersed boundary method is the construction of the height function that gives the distance between the immersed boundaries.

The movies referenced throughout the manuscript are available at author TGF's website at \url{http://scholar.harvard.edu/tfai/home}.

\section{Formulation}
\label{sec:form}

First, we review the standard immersed boundary method, so that the new features of the lubricated method will be clear when they are introduced later on.

\subsection{Classical Immersed Boundary Method}
\label{sec:form_fl}

Consider a domain $\Omega$ filled with incompressible fluid of density $\rho$ and viscosity $\mu$, having velocity $\mb{u}(\mb{x},t)$ and pressure $p(\mb{x},t)$ defined in terms of cartesian coordinates $\mb{x}$. For many relevant problems, the nonlinear term in the Navier-Stokes equations may be neglected, and we will restrict attention to this case in which the fluid is described by the unsteady Stokes equations:
\begin{align}
 \rho\frac{\partial\mb{u}}{\partial t}+\nabla p &= \mu \Delta \mb{u}+\mb{f}\\
 \nabla \cdot \mb{u} &= 0,
\end{align}
where $\mb{f}(\mb{x},t)$ is the applied force density on the fluid.
We are interested in the case in which an immersed elastic body exerts forces on the fluid and thereby determines $\mb{f}$.

The elastic body is described using Lagrangian coordinates $\mb{q}$, e.g.\ arclength in the undeformed frame for a one-dimensional immersed boundary, and the function $\mb{X}(\mb{q},t)$ gives the cartesian position at time $t$ of the material point labeled by $\mb{q}$ along the boundary. We assume that the Lagrangian force density $\mb{F}(\mb{q},t)$ is calculated in terms of $\mb{X}(\mb{q},t)$; for instance, if there exists an energy functional $E = E[\mb{X}(\cdot,t)]$, $\mb{F}$ may be calculated via $\mb{F}=-\varrho E/\varrho \mb{X}$, where $\varrho /\varrho \mb{X}$ denotes the first variation or Fr\'{e}chet derivative. The Dirac delta function is used to convert the Eulerian force density from Lagrangian to cartesian coordinates via
\begin{equation}\label{eqn:spread}
 \mb{f}(\mb{x},t)=\int_\mb{q} \mb{F}(\mb{q},t)\delta(\mb{x}-\mb{X}(\mb{q},t))\ud \mb{q}.
\end{equation}
The immersed boundary moves at a velocity $\mb{U}^\text{IB}(\mb{q},t)$ equal to the local fluid velocity. This condition can be expressed in terms of the Dirac delta function as follows:
\begin{equation}\label{eqn:interp}
 \mb{U}^\text{IB}(\mb{q},t)=\int_\Omega \mb{u}(\mb{x},t)\delta(\mb{x}-\mb{X}(\mb{q},t))\ud \mb{x}.
\end{equation}
As noted in \cite{acta}, the above equation can be simplified using the definition of the Dirac delta function to obtain $\mb{U}^\text{IB}(\mb{q},t)=\mb{u}(\mb{X}(\mb{q},t),t)$, but the expression \eqref{eqn:interp} is particularly convenient since this integral form is the basis of the numerical method used, in which the Dirac delta functions are replaced by regularized delta functions and $\mb{U}^\text{IB}(\mb{q},t)$ is evaluated by summation over the fluid grid. Regularization gives the immersed boundaries an effective width or blurriness on the order of the grid spacing $\Delta x$.

\subsection{Lubrication Corrections}
\label{sec:form_lube}

Now, consider two immersed boundaries separated by a thin layer of fluid with height $h \ll 1$. Note that throughout this article $h$ will be used to represent height and should \emph{not} be assumed equal to the grid spacing $\Delta x$ (or $\Delta y$). In the limit $h \to 0$, the two immersed boundaries coincide and move at identical velocities given by \eqref{eqn:interp}. By the same logic, we expect that when boundary velocities are interpolated using regularized delta functions, boundaries separated by a gap smaller than the grid spacing $\Delta x$ will move at similar velocities and stick together. A thought experiment illustrates the issue: imagine two parallel lines separated by a distance $h$ and pulled in opposite directions. The Lagrangian forces on the lines are equal and opposite, and if the regularization length scale is large compared to the spacing between the lines, the delta function layers of force will be blurred so that they nearly cancel. In the limit $h \to 0$, this leads to a vanishing Eulerian force density in \eqref{eqn:spread}, and therefore vanishing velocity fields $\mb{u}(\mb{x},t)$ and $\mb{U}^\text{IB}(\mb{q},t)$. In practice, similar difficulties are encountered whenever $h \lesssim \Delta x$, since in this case the spread forces partially cancel and the thin lubricating layers separating immersed boundaries are not sufficiently resolved. Instead of computing the distinct velocities of the two boundaries, underresolved simulations assign the averaged velocity to both immersed boundaries, and the lines fail to move unless an extremely large force is exerted.

\begin{figure}[!ht]
         \begin{center}
	 \includegraphics[width=3 in]{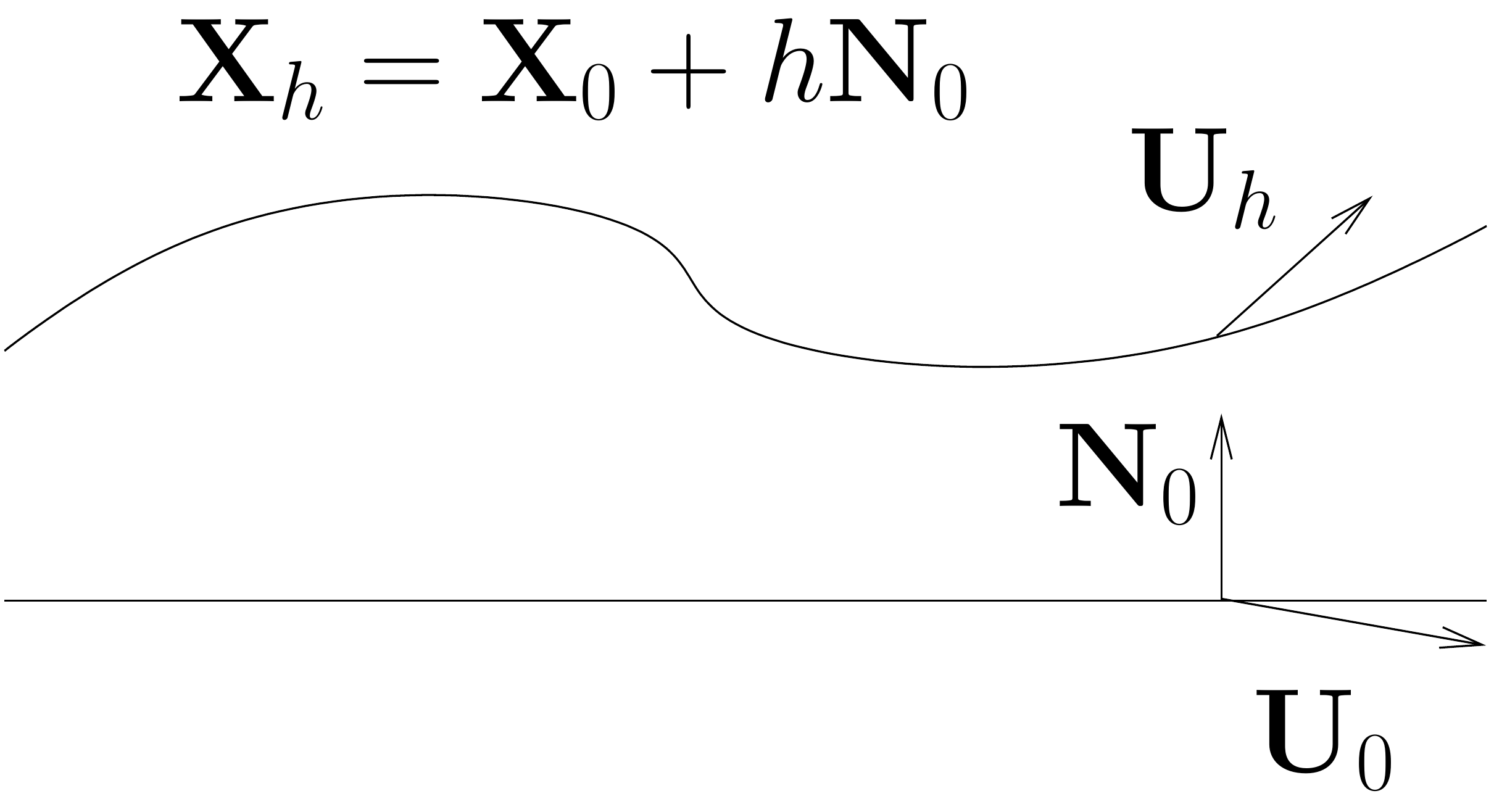}
         \caption{Schematic of the variables used in the lubrication immersed boundary method.}
         \label{fig:test1_schem1}
         \end{center}
\end{figure}
Motivated by this observation, we propose an alternative formulation for the advection velocity $\mb{U}(\mb{q},t)$ with the property $\mb{U}(\mb{q},t) \neq \mb{U}^\text{IB}(\mb{q},t)$. We use the average velocity computed by the standard immersed boundary method, but compute the velocity difference in an alternative way. Assume that there are two immersed boundaries close to touching, and let $\mb{U}_0$ and $\mb{U}_h$ be the velocities on the lower and upper surfaces, respectively, with corresponding force densities $\mb{F}_0$ and $\mb{F}_h$ and positions $\mb{X}_0$ and $\mb{X}_h$ satisfying $\mb{X}_h = \mb{X}_0+h\mb{N}_0$, where $\mb{N}_0$ is the unit normal to the lower surface at position $\mb{X}_0$ (Figure \ref{fig:test1_schem1}). In our formulation, $\mb{U}_0$ and $\mb{U}_h$ are determined by the following equations specifying their sum and difference:
\begin{align}
 &\mb{U}_h+\mb{U}_0 = \mb{U}_h^\text{IB}+\mb{U}_0^\text{IB} \label{eq:lubesum}\\
  ( &\mb{U}_h-\mb{U}_0)\cdot\mb{T}_0 = \frac{h}{2}\Bigg(\frac{(\mb{N}_0\cdot\mb{N}_h)^2\edone{(F_h^\parallel/J_h)}-\edone{(F_0^\parallel/J_0)}}{\mu}\notag\\
  &\qquad\qquad\qquad\qquad\qquad+\edone{\left((\mb{N_0} \cdot \nabla)\mb{u}\right)}\Big |_{0^-}\mkern-18mu\cdot\mb{T}_0+\edone{\left((\mb{N_0} \cdot \nabla)\mb{u}\right)}\Big |_{h^+}\mkern-18mu\cdot\mb{T}_0\Bigg) \label{eq:lubediff}\\
     ( &\mb{U}_h-\mb{U}_0)\cdot\mb{N_0}= \edone{(\mb{T_0} \cdot \nabla)}\left(\frac{h^3}{24\mu} \edone{(\mb{T_0} \cdot \nabla)}\left( \edone{\frac{F_0^\perp}{J_0}}-\edone{\frac{F_h^\perp}{J_h}}+p |_{h^+}+p |_{0^-}\right)\right)\notag\\
       &\qquad \qquad \qquad -\edone{(\mb{T_0} \cdot \nabla)}\left(h\mb{T}_0\cdot\left(\frac{\mb{U}_h+\mb{U}_0}{2}\right)\right)  +\left(\mb{U}_h\cdot\mb{T}_0\right) (\edone{(\mb{T_0} \cdot \nabla)} h), \label{eq:lubediffperp}
\end{align}
where $\mb{T}$ is the unit tangent vector and $J = \det\left(\partial X_i/ \partial q_j\right)$ is the Jacobian, so that $\mb{F}/J$ has units of an Eulerian force density. We write $F^\parallel = \mb{T}\cdot\mb{F}$ and $F^\perp = \mb{N}\cdot\mb{F}$ to denote the tangential and normal components of the force density, respectively. \edone{The expressions $(\mb{N_0} \cdot \nabla)$ and $(\mb{T_0} \cdot \nabla)$ are used to denote directional derivatives normal and tangent to the lower surface, respectively. In the case of a flat lower boundary aligned with the $x$-axis, $(\mb{N_0} \cdot \nabla) \equiv \partial/\partial y$, and $(\mb{T_0} \cdot \nabla) \equiv \partial/\partial x$.}

Together, \eqref{eq:lubesum} and \eqref{eq:lubediff}--\eqref{eq:lubediffperp} define $\mb{U}_h$ and $\mb{U}_0$; to summarize, the standard immersed boundary formulation is used for the average velocity and a sharp formulation is used for the velocity difference. It is reasonable to compute the mean velocity using the standard immersed boundary method as in \eqref{eq:lubesum}, since smoothing tends to yield the net motion as in the thought experiment above. The distinction between our approach and the standard immersed boundary method is the formulation for the difference $\mb{U}_h-\mb{U}_0$, which has several practical advantages compared to the standard method. Consider \eqref{eq:lubediff}, for instance: the Lagrangian forces $\mb{F}_0$ and $\mb{F}_h$ are quantities supported on the immersed boundary, and it follows that $F_0^\parallel$ and $F_h^\parallel$ are never smoothed by convolution against a delta function kernel. Moreover, the quantities \edone{$((\mb{N_0} \cdot \nabla)\mb{u})|_{0^-}$} and \edone{$\left((\mb{N_0} \cdot \nabla\right)\mb{u}) |_{h^+}$} are \emph{outer} derivatives, which implies they are typically well-resolved and/or small in magnitude. Therefore, the terms on the right-hand side of \eqref{eq:lubediff} are sharply defined and convenient to use, and \eqref{eq:lubediffperp} is similarly formulated in terms of quantities supported on the boundary or outside the lubrication layer. Next, we show how \eqref{eq:lubediff}--\eqref{eq:lubediffperp} may be derived from lubrication theory. 

\subsection{Lubrication Equations}
\label{sec:deriv}

In the thought experiment presented above involving two parallel lines pulled in opposite directions with $x$-velocities $U_0$ and $U_h$, the exact solution $\mb{u}=(u,v)$ inside the gap between the lines is linear shear flow, with $u(y) = (U_h-U_0) y/h+U_0$ and $v \equiv 0$. We wish to use this formula for $u(y)$ to obtain the difference $U_h-U_0$ in terms of known quantities, and it will turn out that the derivative $\partial u/\partial y$ normal to the lines is a convenient choice. By the fundamental theorem of calculus,
\begin{equation}
 U_h-U_0 = \int_0^h \frac{\partial u}{\partial y} \ud y,
\end{equation}
and because the derivative $\partial u/\partial y$ is constant for linear shear, it follows that 
\begin{equation}
\label{eq:diff_ss}
 U_h-U_0=h \frac{\partial u}{\partial y} \Bigg |_{y^*},
\end{equation}
where $y^*$ is any value of $y$ in $(0,h)$. Although the above equation is trivial, as it simply states that the flow in a linear shear is completely determined by its shear rate, surprisingly there is a similar statement for flows through narrow lubrication layers in general.

To prove this more general statement, we first review some classical results from lubrication theory, although we will use these results in a nonstandard manner. In most applications of lubrication theory, velocity boundary conditions are prescribed and the pressure is calculated by solving the Reynolds equation \cite{Acheson1990}, whereas in our case the relative velocity of the boundaries is the unknown.

Assume a flat lower surface aligned with the $x$-direction is separated by height $h$ from a (possibly curved) upper surface with characteristic horizontal length scale $L$, where $h \ll L$. The unsteady Stokes equations can be simplified substantially in this case by neglecting terms of order $h/L$ and smaller to obtain \cite{Acheson1990}
\begin{align}
\frac{\partial p}{\partial x} &= \mu \frac{\partial^2 u}{\partial y^2} \label{eq:lube1}\\
\frac{\partial p}{\partial y} &=0\\
\partial u/\partial x+\partial v/\partial y &= 0.
\end{align}
Since $p$ is $y$-independent, \eqref{eq:lube1} may be integrated explicitly to yield
\begin{align}
\label{eq:dudy_lube}
 \frac{\partial u}{\partial y} = \frac{1}{2\mu}\frac{\partial p}{\partial x} (2y-h) + \frac{U_h-U_0}{h}.
\end{align}
Because the first term of the right-hand side has opposite signs for $y=0$ and $y=h$, we may eliminate the pressure gradient and invert the above equation for $U_h-U_0$ to get
\begin{align}
\label{eq:diff_lu1}
  U_h-U_0=\frac{h}{2}\left(\frac{\partial u}{\partial y}\Big|_{0^+}+\frac{\partial u}{\partial y}\Big|_{h^-}\right).
\end{align}
This implies that the observation \eqref{eq:diff_ss} made in the special case of linear shear that the velocity difference may be computed by knowledge of $\partial u/\partial y$ at the boundaries applies to lubrication flows in general. In fact, the pressure through the gap can also be recovered using
\begin{align}
\label{eq:diff_lu2}
  \frac{\partial p}{\partial x} = \frac{\mu}{h}\left(\frac{\partial u}{\partial y}\Big|_{h^-}-\frac{\partial u}{\partial y}\Big|_{0^+}\right).
\end{align}
Knowledge of the pressure gradient and boundary velocities is sufficient to reconstruct the full tangential velocity profile through the gap by the formula
\begin{align}
 \label{eq:u_lube}
u(y) = \frac{1}{2\mu}\frac{\partial p}{\partial x}y(y-h)+\frac{U_h-U_0}{h}y+U_0,
\end{align}
which follows from integrating \eqref{eq:dudy_lube}. (The normal velocities may be computed in terms of the pressure gradient and tangential velocities, as shown later on.)

We have now derived how the boundary velocities and flow through the gap may be calculated in terms of normal derivatives on the boundaries. Next, we will show these normal derivatives may themselves be calculated in terms of more convenient quantities.

\subsection{Jump Conditions}
\label{sec:jump}

Here, we review the classical jump conditions across an immersed boundary \cite{lai2001remark,xu2006systematic}, closely following the exposition of \cite{lecjump}. Let $[\phi]$ denote the jump in a quantity $\phi$ across an immersed boundary, so that e.g.~$[\sigma_{ij}]$ for $i,\,j = 1,\,2$ gives the jump in the components of the stress tensor. Since immersed boundaries are considered massless, all stresses must balance, i.e.~
\begin{equation}
 [\sigma_{ij}]N_j+F_i/J=0,
\end{equation}
where $\mb{F}$ is the Lagrangian force density, $\mb{N}$ is the unit normal vector, and $J$ is Jacobian. \edone{We will continue to assume for now that the lower surface is aligned with the $x$-axis so that $\mb{N}_0 \equiv \mathbf{\hat{\j}}$. However, the upper surface is not assumed flat, i.e.~$\mb{N}_h$ may vary from point to point depending on its curvature. As we will need to use the jump conditions across both surfaces, we use a more general notation that involves surface normals and tangents.} In a Newtonian fluid the stress tensor is given by
\begin{equation}
\sigma_{ij}=-p\delta_{ij}+\mu\left(\partial u_i/\partial x_j+\partial u_j/\partial x_i \right).
\end{equation}
It follows that
\begin{equation}
-[p]N_i+\mu \left[ \partial u_i/\partial x_j \right]N_j +\mu \left[ \partial u_j/\partial x_i \right]N_j +F_j/J=0,
\end{equation}
 which can be simplified \cite{lecjump} to
\begin{equation}
\label{eq:jump}
[p]\mb{N}-\mu \left[  (\mb{N} \cdot \nabla) \mb{u}\right] = \mb{F}/J.
\end{equation}
Separating \eqref{eq:jump} into tangential and normal components, we have
\begin{align}
-\mu \left[ (\mb{N} \cdot \nabla) \mb{u}\right] &= \left(\mb{F}-(\mb{F}\cdot\mb{N})\mb{N}\right)/J
 \label{eq:jump1}\\
[p]&= \left(\mb{F}\cdot\mb{N}\right)/J, \label{eq:jump2}
\end{align}
where we have used that the normal component of $\left[  (\mb{N} \cdot \nabla) \mb{u}\right]$ vanishes \cite{lecjump}. The first jump condition \eqref{eq:jump1} relates forces on the immersed boundary to the normal derivative in velocity, which by \eqref{eq:diff_lu1} is just the quantity needed to find $U_h-U_0$. \edone{On the lower surface, $\mb{N}_0 \equiv \mathbf{\hat{\j}}$ implies that $(\mb{N}_0 \cdot \nabla) \mb{u} = \partial \mb{u}/\partial y$, and we may use \eqref{eq:jump1} to rewrite the first term of \eqref{eq:diff_lu1} as
\begin{align}
 \frac{\partial u}{\partial y}\Big|_{0^+} &= -\frac{F_0^\parallel/J_0}{\mu}+\frac{\partial u}{\partial y}\Bigg |_{0^-},
  \label{eq:jump_low}
\end{align}
where $F^\parallel = (\mb{F}\cdot\mb{T})$ is the tangential component of the Lagrangian force density and $\mb{T}$ is a unit tangent vector, i.e.~$\mb{T}_0 \equiv \mathbf{\hat{\i}}$. On the other hand, the unit normal $\mb{N}_h$ to the upper surface does not in general point in the $y$-direction, and so the associated jump condition does not directly involve $\partial \mb{u}/\partial y$. However, since 
\begin{align}
 \mb{N}_0 = (\mb{N}_0\cdot\mb{N}_h)\mb{N}_h+(\mb{N}_0\cdot\mb{T}_h)\mb{T}_h,
\end{align}
we have for the jump across the upper surface that
\begin{align}
  \left[ \partial \mb{u}/\partial y\right] = \left[ (\mb{N_0} \cdot \nabla) \mb{u}\right] &= (\mb{N}_0\cdot\mb{N}_h)\left[ (\mb{N_h} \cdot \nabla) \mb{u}\right]+(\mb{N}_0\cdot\mb{T}_h)\left[ (\mb{T_h} \cdot \nabla) \mb{u}\right] \notag\\
 &= (\mb{N}_0\cdot\mb{N}_h)\left[ (\mb{N_h} \cdot \nabla) \mb{u}\right] \notag\\
 &=(\mb{N}_0\cdot\mb{N}_h)\left(\mb{F}_h-(\mb{F}_h\cdot\mb{N}_h)\mb{N}_h\right)/J_h
\end{align}
where we have used that there is no jump in tangential derivatives \cite{lecjump} in the second line, and the jump condition \eqref{eq:jump1} in the final line. The right hand side above is clearly tangential to the upper surface, and therefore we may take the dot product with $\mb{T}_0$ to obtain $[\partial u/\partial y]$. Since $(\mb{T}_0\cdot\mb{T}_h)=(\mb{N}_0\cdot\mb{N}_h)$ in 2D this results in
\begin{align}
  \frac{\partial u}{\partial y}\Big|_{h^-} &= (\mb{N}_0\cdot\mb{N}_h)^2\frac{F_h^\parallel/J_h}{\mu}+\frac{\partial u}{\partial y}\Bigg |_{h^+}.
  \label{eq:jump_up}
\end{align}
}
Substituting into \eqref{eq:diff_lu1} the expressions \eqref{eq:jump_low} and \eqref{eq:jump_up} yields
\begin{equation}
\label{eq:lubediff_simple}
   U_h-U_0 = \frac{h}{2}\left(\frac{(\mb{N}_0\cdot\mb{N}_h)^2\edone{(F_h^\parallel/J_h)}-\edone{(F_0^\parallel/J_0)}}{\mu}+\frac{\partial u}{\partial y}\Bigg |_{0^-}+\frac{\partial u}{\partial y}\Bigg |_{h^+}\right),
\end{equation}

To this point, we have not considered the normal velocities $V_0$ and $V_h$ on the lower and upper boundaries, respectively. In order to obtain an expression for $V_h-V_0$, we invoke the incompressibility constraint $\partial u/\partial x+\partial v/\partial y=0$ to write
\begin{align}
 V_h-V_0 &= \int_0^h \frac{\partial v}{\partial y} \ud y \notag\\
  &= -\int_0^h \frac{\partial u}{\partial x} \ud y \notag \\
   &= -\left(\frac{\partial q}{\partial x}-U_h \frac{\partial h}{\partial x} \right) \label{eq:diff_n1},
\end{align}
where $q = \int_0^h u \ud y$ is the flow rate and we have used the fundamental theorem of calculus in the first and third equalities. Integrating \eqref{eq:u_lube} directly yields
\begin{equation}
 q = \frac{-h^3}{12 \mu}\frac{\partial p}{\partial x}+h\frac{U_h+U_0}{2},
\end{equation}
and substituting this expression into \eqref{eq:diff_n1} results in
\begin{equation}
 V_h-V_0 = \left(\frac{h^3}{12\mu} \frac{\partial p}{\partial x}\right)_x-\left(h\left(\frac{U_h+U_0}{2}\right)\right)_x+U_h \frac{\partial h}{\partial x}.
 \label{eq:u_lube2_int}
\end{equation}
\edone{Note the asymmetry with respect to the lower and upper surfaces, in that the final term on the right-hand side of \eqref{eq:u_lube2_int} has no corresponding term involving $U_0$. This is a consequence of defining the height in terms of the lower surface normal.}

As was the case for the tangential velocities, we would like to be able to compute $V_h-V_0$ in terms of convenient quantities such as the boundary forces. This is possible using the above formula \eqref{eq:u_lube2_int}, which involves tangential velocities already known through \eqref{eq:lubediff_simple}, the pressure gradient $\partial p/\partial x$ already known through \eqref{eq:diff_lu2}, and derivatives of these quantities. In practice, to compute the pressure gradient we have found it preferable to use the jump condition \eqref{eq:jump2}, according to which
\begin{align}
 p |_{0^+} &= \edone{\frac{F_0^\perp}{J_0}}+p |_{0^-}\\
 p |_{h^-} &= -\edone{\frac{F_h^\perp}{J_h}}+p |_{h^+},
\end{align}
where $F^\perp = (\mb{F}\cdot\mb{N})$. We average $p |_{0^+}$ and $p |_{h^-}$ to get a pressure in the gap and obtain the following formula:
\begin{equation}
 V_h-V_0 = \left(\frac{h^3}{24\mu} \left( \edone{\frac{F_0^\perp}{J_0}}-\edone{\frac{F_h^\perp}{J_h}}+p |_{h^+}+p |_{0^-}\right)_x\right)_x-\left(h\left(\frac{U_h+U_0}{2}\right)\right)_x+U_h \frac{\partial h}{\partial x}.
   \label{eq:u_lube2}
\end{equation}
Thus we have an equation for $V_h-V_0$ that involves boundary forces, horizontal velocities that are themselves computed in terms of boundary forces, and the outer pressures $p |_{0^-}$ and $p |_{h^+}$.

\edone{For the general case of a nonlinear lower surface}, we replace $x$- and $y$-velocities by velocities tangential and normal to the lower surface, respectively, and $x$-derivatives by tangential derivatives to the lower surface, resulting in the equations \eqref{eq:lubediff} and \eqref{eq:lubediffperp} for the lubricated immersed boundary method. \edone{Substituting coordinates in this manner is reasonable provided that $h \ll L$ so that the effects of curvature are negligible.} The equations analogous to \eqref{eq:lubediff_simple} and \eqref{eq:u_lube2} in the general case are \eqref{eq:lubediff} and \eqref{eq:lubediffperp} of Section \ref{sec:form_lube}, respectively. The equation analogous to \eqref{eq:u_lube} for the tangential velocity through the gap in the case of curved surfaces is
\begin{align}
\mb{u}(\mb{x})\cdot\mb{T}_0 = \frac{1}{2\mu} w(w-h) (\mb{T}_0 \cdot \nabla)p+w\left(\frac{\mb{U}_h-\mb{U}_0}{h}\right)\cdot\mb{T}_0+\mb{U}_0\cdot\mb{T}_0,  \label{eq:u_lube3}
\end{align}
where $w \in [0,h]$ and where $\mb{x}=\mb{X}_0+w\mb{N}_0$ for fixed $\mb{X}_0$, so that $\mb{x}$ travels through the gap in the direction normal to the lower curve. Given this explicit representation for the tangential velocity, the normal component of the gap velocity may be computed as before via incompressibility:
\begin{align}
\mb{u}(\mb{x})\cdot\mb{N}_0 &= \int_0^{\left(\mb{x}\cdot\mb{N}_0\right)} \frac{\partial\left(\mb{u}\cdot\mb{N}_0\right)}{\partial \mb{N}_0} \ud w+\mb{U}_0\cdot\mb{N}_0 \notag \\
&= - \int_0^{\left(\mb{x}\cdot\mb{N}_0\right)} \frac{\partial\left(\mb{u}\cdot\mb{T}_0\right)}{\partial \mb{T}_0} \ud w+\mb{U}_0\cdot\mb{N}_0, \label{eq:u_lube4}
\end{align}
where the integrand can be computed through \eqref{eq:u_lube3}. \edone{Note that for the upper limit of integration in \eqref{eq:u_lube4}, $\mb{x}\cdot\mb{N}_0$ is simply the height above the lower surface of an arbitrary point $\mb{x}$ inside the gap.}



\section{Numerical Method}

\subsection{Spatial discretization}
We next specialize to the case of a rectangular unit domain $\Omega$ of size $L_x$ by $L_y$ and periodic boundary conditions in both directions. We discretize $\Omega$ by overlaying an $N_x \times N_y$ grid with grid spacing $\Delta x=L_x/N_x$ and $\Delta=L_y/N_y$ in the $x$- and $y$-directions, respectively. The $(i,j)$\textsuperscript{th} grid cell has its center at $\mb{x}_{i,j} = ((i+1/2)\Delta x, (j+1/2)\Delta y)$ for $i=0,\dots,N_x-1$ and $j=0,\dots,N_y-1$.

We will consider a one-dimensional boundary immersed in $\Omega$ and parameterized by the Lagrangian coordinate $q \in [0, L_q]$. The immersed boundary is discretized using $N_q$ points so that the Lagrangian grid spacing is $\Delta q = L_q/N_q$, and the positions of Lagrangian points in cartesian coordinates are given by $\mb{X}_k(t) := \mb{X}(k\Delta q, t)$ for $k=0,\dots,N_q-1$.

The grid velocities $\mb{u}_{i,j}$ are located at cell centers $\mb{x}_{i,j}$, as are the external forces. The pressures are located at cell corners, with $p_{i,j}$ defined at the bottom-right corner of the $(i,j)$\textsuperscript{th} grid cell. \edone{Our discretization is based on \cite{almgren1996numerical}, which uses a finite element projection to obtain a symmetric positive definite system that can be solved efficiently by methods such as geometric multigrid.} The unsteady Stokes equations are discretized in space in a centered, 2\textsuperscript{nd}-order accurate manner:
\begin{align}
\rho \frac{\ud \mb{u}_{i,j}}{\ud t} + (\mb{G}p)_{i,j} &= \mu (L\mb{u})_{i,j} +\mb{f}_{i,j}\\
\left(\mb{D} \cdot \mb{u}\right)_{i,j} &= 0,
\end{align}
where
\begin{align}\
 (\mb{G}p)_{i,j} &= \oh \parenVectorstack{{(p_{i+1,j}-p_{i,j}+p_{i+1,j+1}-p_{i,j+1})/\Delta x} {(p_{i,j+1}-p_{i,j}+p_{i+1,j+1}-p_{i+1,j})/\Delta y}}, \\
 (L\mb{u})_{i,j} &= \parenVectorstack{{(u_{i-1,j}-2u_{i,j}+u_{i+1,j})/(\Delta x)^2} {(v_{i,j-1}-2v_{i,j}+v_{i,j+1})/(\Delta y)^2}},\\
 \left(\mb{D} \cdot \mb{u}\right)_{i,j} &=\oh\left(\frac{u_{i,j}-u_{i-1,j}+u_{i,j-1}-u_{i-1,j-1}}{\Delta x}\right.\notag\\
 &\qquad \left. +\frac{v_{i,j}-v_{i,j-1}+v_{i-1,j}-v_{i-1,j-1}}{\Delta y}\right),
\end{align}
and where all index arithmetic is to be interpreted modulo $N_x$ or $N_y$. We also have the equation of motion for the immersed boundary points
\begin{align}
\frac{\ud \mb{X}_{k}}{\ud t} = \mb{U}_k.
\end{align}
To close the system, we must provide discretized equations for the spreading and interpolation operators used to compute $\mb{f}_{i,j}$ and $\mb{U}_k$. They are as follows:
\begin{align}
 \mb{f}_{i,j}(t) &= \sum_{k=0}^{N_q-1} \mb{F}_k(t) \delta_h\left(\mb{x}_{i,j}-\mb{X}_k(t)\right)\Delta q \label{eq:spread_discrete}\\
  \mb{U}_k(t) &= \sum_{i=0}^{N_x-1}\sum_{j=0}^{N_y-1} \mb{u}_{i,j}(t) \delta_h\left(\mb{x}_{i,j}-\mb{X}_k(t)\right)\Delta x\Delta y, \label{eq:interp_discrete}
\end{align}
where we use the standard 4-point immersed boundary delta function $\delta_h(\mb{x}) = \phi(x/\Delta x)\phi(y/\Delta y)$ \cite{acta}, in which $\phi$ is given by
\begin{align}
\phi(r)&=
\begin{cases}
 \frac{1}{8}\left( 5+2r-\sqrt{-7-12r-4r^2} \right), & -2 \le r < -1, \\
\frac{1}{8}\left( 3+2r+\sqrt{1-4r-4r^2} \right), & -1 \le r < 0, \\
\frac{1}{8}\left( 3-2r+\sqrt{1+4r-4r^2} \right), & 0 \le r < 1, \\
\frac{1}{8}\left( 5-2r-\sqrt{-7+12r-4r^2} \right), & 1 \le r \le 2, \\
0, & |r| > 2.
\end{cases}
\label{eq:phi}
\end{align}
We assume that the Lagrangian force $\mb{F}_k(t)$ can be calculated in terms of the positions $\mb{X}_k(t)$, specific examples of which will be given later on.

\subsection{Timestepping}
Upon discretizing in space, we have a semi-discrete problem with unknowns $\mb{u}_{i,j}(t)$, $p_{i,j}(t)$, and $\mb{X}_k(t)$ that depend on a continuous time variable. The problem is fully discretized by introducing a timestep $\Delta t$ and letting $\mb{u}^n_{i,j} := \mb{u}_{i,j}(n\Delta t)$, $p^n_{i,j} := p_{i,j}(n\Delta t)$, and $\mb{X}^n_k := \mb{X}_k(n\Delta t)$. We use a second-order accurate predictor-corrector timestepping scheme. A preliminary approximation $\mb{X}^*$ to $\mb{X}^{n+1}$ is first computed through a first-order accurate forward Euler step:
\begin{equation}
 \mb{X}^* = \mb{X}^n+\Delta t \mb{U}^n,
\end{equation}
where $\mb{U}^n$ is computed in terms of the grid velocity $\mb{u}^n$ via \eqref{eq:interp_discrete}.
We calculate $\mb{u}^{n+1}$ and $p^{n+1/2}$ by solving a time-centered system of equations:
\begin{align}
\rho \frac{\mb{u}^{n+1}-\mb{u}^n}{\Delta t} + \mb{G}p^{n+1/2} &= \frac{\mu}{2} L\left(\mb{u}^n+\mb{u}^{n+1}\right) +\oh\left(\mb{f}^n+\mb{f}^*\right) \label{eq:us_corr_mom}\\
\mb{D} \cdot \mb{u}^{n+1} &= 0, \label{eq:us_corr_div}
\end{align}
where $\mb{f}^*$ is calculated from $\mb{X}^*$. Finally, we update the positions in a time-centered manner to obtain $\mb{X}^{n+1}$:
\begin{equation}
  \mb{X}^{n+1} = \mb{X}^n+\frac{\Delta t}{2}\left( \mb{U}^n+\mb{U}^*\right),
\end{equation}
where $\mb{U}^*$ is computed in terms of $\mb{X}^*$. Although this scheme involves the first-order accurate preliminary guess $\mb{X}^*$, the approach is similar to the broader class of Runge-Kutta schemes and it can be shown to be second-order accurate (see for instance the related scheme in \cite{vc_ib,vc_ib2}).

To actually solve for the velocity and pressure in \eqref{eq:us_corr_mom}, we take the discrete divergence of both sides of the equation and use \eqref{eq:us_corr_div} to eliminate the discretely divergence-free velocity field (we assume the initial velocity satisfies $\mb{D} \cdot \mb{u}^0 = 0$). This results in a Poisson problem for the pressure $p^{n+1/2}$, and once this problem has been solved \eqref{eq:us_corr_mom} may be rewritten as an explicit equation for $\mb{u}^{n+1}$. We solve this linear system using geometric multigrid \cite{brandt1977multi,briggs2000multigrid}. 

\subsection{Computing the height function}
\label{sec:height}

Given the sensitivity of the underlying equations to the height, e.g.~the dependence in \eqref{eq:u_lube2} of $V_h-V_0$ on $h^3$, the accurate computation of the height function is essential to this method. This additional computational cost compared to the standard immersed boundary method, in which computing the height is not necessary, is the price paid for resolving velocity gradients across thin fluid layers.

\begin{figure}[!ht]
         \begin{center}
	 \includegraphics[width=3 in]{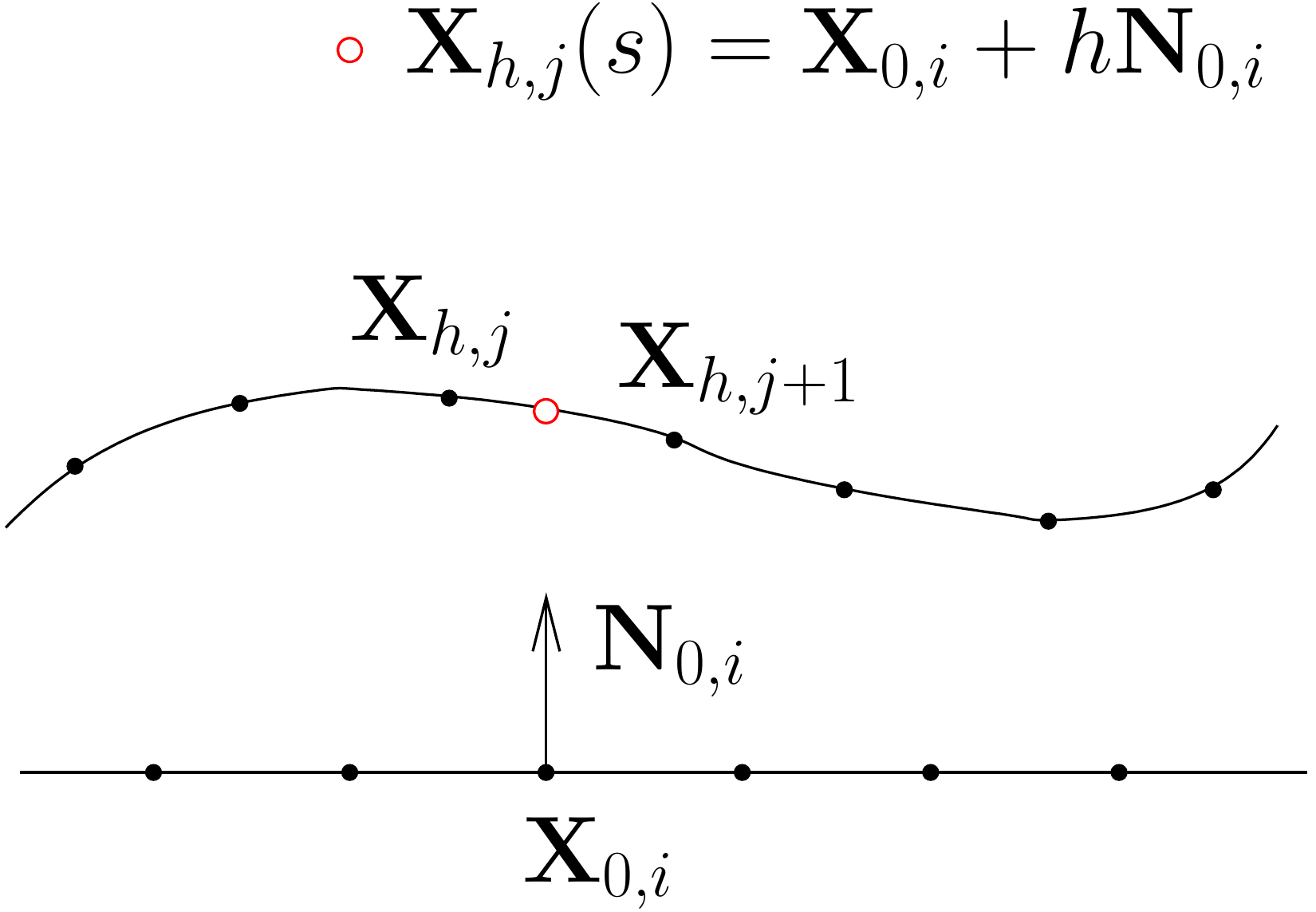}
         \caption{\edone{Computing the height $h$ above a point $\mb{X}_{0,i}$ on the discrete lower surface. The height is defined with respect to the lower surface normal $\mb{N}_{0,i}$. The corresponding point on the upper surface is denoted by $\mb{X}_{h,j}(s)$, where $\mb{X}_{h,j}(s) = \mb{a}_0+\mb{a}_1 s+\mb{a}_2 s^2 +\mb{a}_3 s^3$ for $s \in [0,1]$. That is, the point $\mb{X}_{h,j}(s)$ is located on a section of the upper surface defined by a cubic curve between the discrete points $\mb{X}_{h,j}$ and $\mb{X}_{h,j+1}$.}}
         \label{fig:height1}
         \end{center}
\end{figure}
To compute the height $h_i$ above a point on the lower surface with index $i$ and position $\mb{X}_{0,i}$ (Figure \ref{fig:height1}), we approximate the upper surface by a piecewise cubic curve. On an arbitrary section of the upper surface with endpoints $\mb{X}_{h,j}$ and $\mb{X}_{h,j+1}$, the curve is given by $\mb{X}_{h,j}(s) = \mb{a}_0+\mb{a}_1 s+\mb{a}_2 s^2 +\mb{a}_3 s^3$ for $s \in [0,1]$. The four unknown constants $\left\{\mb{a}_k \big| k=0,\dots,3\right\}$ are determined by the conditions
\begin{align}
 \mb{X}_h(0) &= \mb{X}_{h,j} \label{eq:cub1}\\
  \mb{X}_h(1) &= \mb{X}_{h,j+1}\\
 \frac{\ud \mb{X}_h}{\ud s}(0) &=  \mb{T}_{h,j}\\
  \frac{\ud \mb{X}_h}{\ud s}(1) &=  \mb{T}_{h,j+1} \label{eq:cub4},
\end{align}
where $\mb{T}_{h,j}$, the tangent vector to the upper surface at point $j$, is defined as the weighted average of the two edge vectors connected to vertex $\mb{X}_{h,j}$:
\begin{equation}
 \mb{T}_{h,j} :=  \frac{\left(\mb{X}_{h,j+1}-\mb{X}_{h,j}\right)\norm{\mb{X}_{h,j+1}-\mb{X}_{h,j}}+\left(\mb{X}_{h,j}-\mb{X}_{h,j-1}\right)\norm{\mb{X}_{h,j}-\mb{X}_{h,j-1}}}{\norm{\mb{X}_{h,j+1}-\mb{X}_{h,j}}+\norm{\mb{X}_{h,j}-\mb{X}_{h,j-1}}}
 \end{equation}
Solving the system \eqref{eq:cub1}--\eqref{eq:cub4} results in the coefficients
\begin{align}
 \mb{a}_0 &= \mb{X}_{h,j}\\
  \mb{a}_1 &= \mb{T}_{h,j}\\
   \mb{a}_2 &=  3\left(\mb{X}_{h,j+1}-\mb{X}_{h,j}\right)-2\mb{T}_{h,j}-\mb{T}_{h,j+1}\\
 \mb{a}_3 &=  \mb{T}_{h,j}+\mb{T}_{h,j+1}-2\left(\mb{X}_{h,j+1}-\mb{X}_{h,j}\right).
\end{align}

Defining components of the tangent vector by $\mb{T} = (T^1,T^2)$, the unit normal $\mb{N}$ is given by
\begin{equation}
 \mb{N} =  (-T^2,T^1)/\norm{\mb{T}}.
 \end{equation}
To calculate the height above $\mb{X}_{0,i}$, we find the pair $(s,h)$ and the index $j$ such that
\begin{equation}
  \mb{X}_{0,i}+h \mb{N}_{0,i} = \mb{X}_{h,j}(s).
  \label{eq:lower_height}
\end{equation}
On a section $j$ of the piecewise cubic upper curve, we find $s$ by taking the dot product of both sides with $\mb{T}_{0,i}$ to obtain
\begin{align}
  \mb{X}_{0,i}\cdot \mb{T}_{0,i} &= \mb{X}_{h,j}(s)\cdot \mb{T}_{0,i} \notag \\
&= \left(\mb{a}_0+\mb{a}_1 s+\mb{a}_2 s^2 +\mb{a}_3 s^3\right)\cdot \mb{T}_{0,i},
\end{align}
and the resulting nonlinear equation for $s$ is solved using Newton's method. Having found $s$, the height $h$ is given by
\begin{equation}
  h = \left(\mb{X}_{h,j}(s)-\mb{X}_{0,i}\right) \cdot \mb{N}_{0,i}.
  \end{equation}
To take advantage of the quadratic convergence of Newton's method, it is crucial to have a good initial guess. We use a piecewise linear (rather than piecewise cubic) approximation to the upper surface to obtain the initial guess $s_0$. This has the virtue of being easy to solve: finding $s_0$ such that $\mb{X}_{0,i}+h_0 \mb{N}_{0,i} = \mb{X}_{h,j}+s_0(\mb{X}_{h,j+1}-\mb{X}_{h,j})$ yields
\begin{equation}
s_0 = \frac{\left(\mb{X}_{0,i}-\mb{X}_{h,j}\right)\cdot\mb{T}_{0,i}}{\left(\mb{X}_{h,j+1}-\mb{X}_{h,j}\right)\cdot\mb{T}_{0,i}}.
\end{equation}

The above procedure is done for each $j$ corresponding to a section of the upper surface. For some applications, there is a unique section that yields $s \in [0,1]$ and $h > 0$. However, in some cases there are multiple sections with positive heights above $\mb{X}_{0,i}$ and $s$ in the appropriate range, such as when the upper surface is a closed curve. In that case, the height is defined to be the minimum of the eligible values. We note that a na\"{i}ve implementation of this algorithm to compute heights would involve a nested loop over $i$ and $j$ and result in an $\mathcal{O}(N^2)$ operation. However, looping over all $j$ can be avoided since the lubricated method is only applicable when two immersed boundaries are separated by a distance comparable to the grid spacing. Only nearby points need to be checked, and this allows for fast methods based on binning points into containers of size on the order of $\Delta x$ by $\Delta y$.

\edone{We apply the lubrication equations \eqref{eq:lubesum}--\eqref{eq:lubediffperp} at a point $\mb{X}_{0,i}$ together with its corresponding point $\mb{X}_{h,j}(s)$ on the upper surface. Since in general $\mb{X}_{h,j}(s)$ does not coincide with a discrete node, this involves defining a velocity $\mb{U}^\text{IB}_{h,j}(s)$ and Lagrangian force $\mb{F}_{h,j}(s)$ at the position $\mb{X}_{h,j}(s)$. We interpolate these quantities from the discrete points on the upper surface by the same cubic scheme used above to compute $\mb{X}_{h,j}(s)$. Note that although we use the lubrication corrections to solve for both $\mb{U}_{0,i}$  and $\mb{U}_{h,j}(s)$, only $\mb{U}_{0,i}$ is actually used. The velocities on the upper surface are calculated in an analogous but separate step.}

We also need to define heights $h_j$ below each point on the upper surface in order to compute velocities there. The procedure used to do so is similar to the one just described for the lower surface, although since the height is defined relative to the lower surface normal, the equation that needs to be solved is slightly more complicated. We describe this procedure in detail in Appendix \ref{app:height}.

\edtwo{
\subsection{Outer derivatives}
To compute outer derivatives of velocity, we use the global velocity field $\mb{u}^\text{IB}(\mb{X},t)$ and the fact that the discrete delta function $\delta_h(\mb{x})$ is continuously differentiable. The global velocity gradient is defined by
\begin{align}
\nabla\mb{u}^\text{IB}(\mb{X},t) &= -\sum_{i=0}^{N_x-1}\sum_{j=0}^{N_y-1} \mb{u}_{i,j}(t) \nabla\delta_h\left(\mb{x}_{i,j}-\mb{X}\right)\Delta x\Delta y,
\end{align}
where $\nabla \delta_h(\mb{x}) = \left(\frac{1}{\Delta x}\phi'(x/\Delta x)\phi(y/\Delta y),\frac{1}{\Delta y} \phi(x/\Delta x)\phi'(y/\Delta y)\right)$ and $\phi'(r)$ is calculated by analytically differentiating the expressions in \eqref{eq:phi} to obtain
\begin{align}
\phi'(r)&=
\begin{cases}
 \frac{1}{8}\left( 2+\frac{6+4r}{\sqrt{-7-12r-4r^2}} \right), & -2 \le r < -1, \\
\frac{1}{8}\left( 2+\frac{-2-4r}{\sqrt{1-4r-4r^2}} \right), & -1 \le r < 0, \\
\frac{1}{8}\left( -2+\frac{2-4r}{\sqrt{1+4r-4r^2}} \right), & 0 \le r < 1, \\
\frac{1}{8}\left( -2+\frac{-6+4r}{\sqrt{-7+12r-4r^2}} \right), & 1 \le r \le 2, \\
0, & |r| > 2.
\end{cases}
\end{align}
As shown in the numerical experiments of Section \ref{sec:simp_shear}, the interpolated gradient gives a good approximation to the outer velocity derivative when it is evaluated at least 2 gridpoints from the boundary. We therefore calculate the outer derivative at a distance of 2 -- 4 gridpoints from the immersed boundaries. The optimal choice depends on the application; for instance, a smaller distance is sometimes needed to avoid evaluating $\nabla \mb{u}^\text{IB}$ outside the domain.
}

\section{Convergence Results}

\subsection{Simple shear}
\label{sec:simp_shear}
As a first demonstration of the lubricated immersed boundary method, we simulate the situation described in Section \ref{sec:form_lube} of two parallel lines separated by a distance $h$ and pulled in opposite directions.
\begin{figure}[!ht]
         \begin{center}
         	\subfloat[Setup]{
	 \includegraphics[width=3 in]{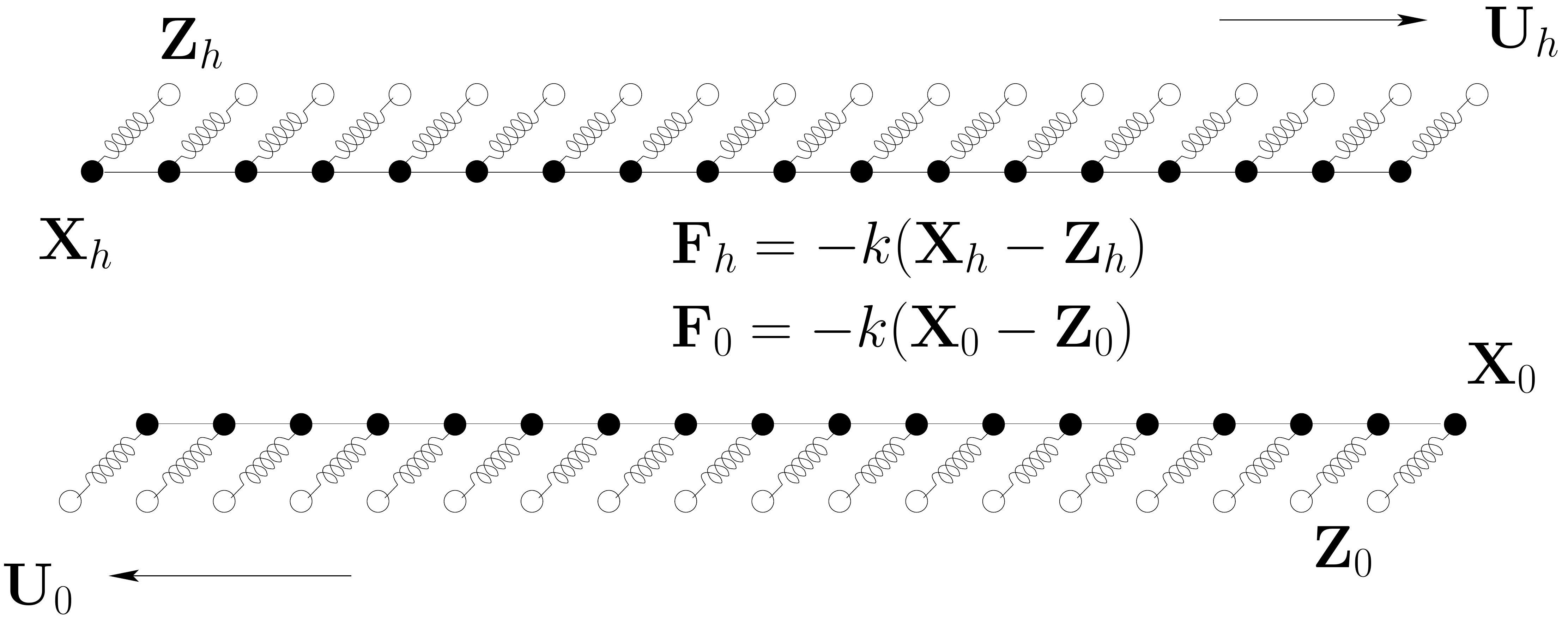}}\\
	 	\subfloat[Velocity profile]{
	 \includegraphics[width=3 in]{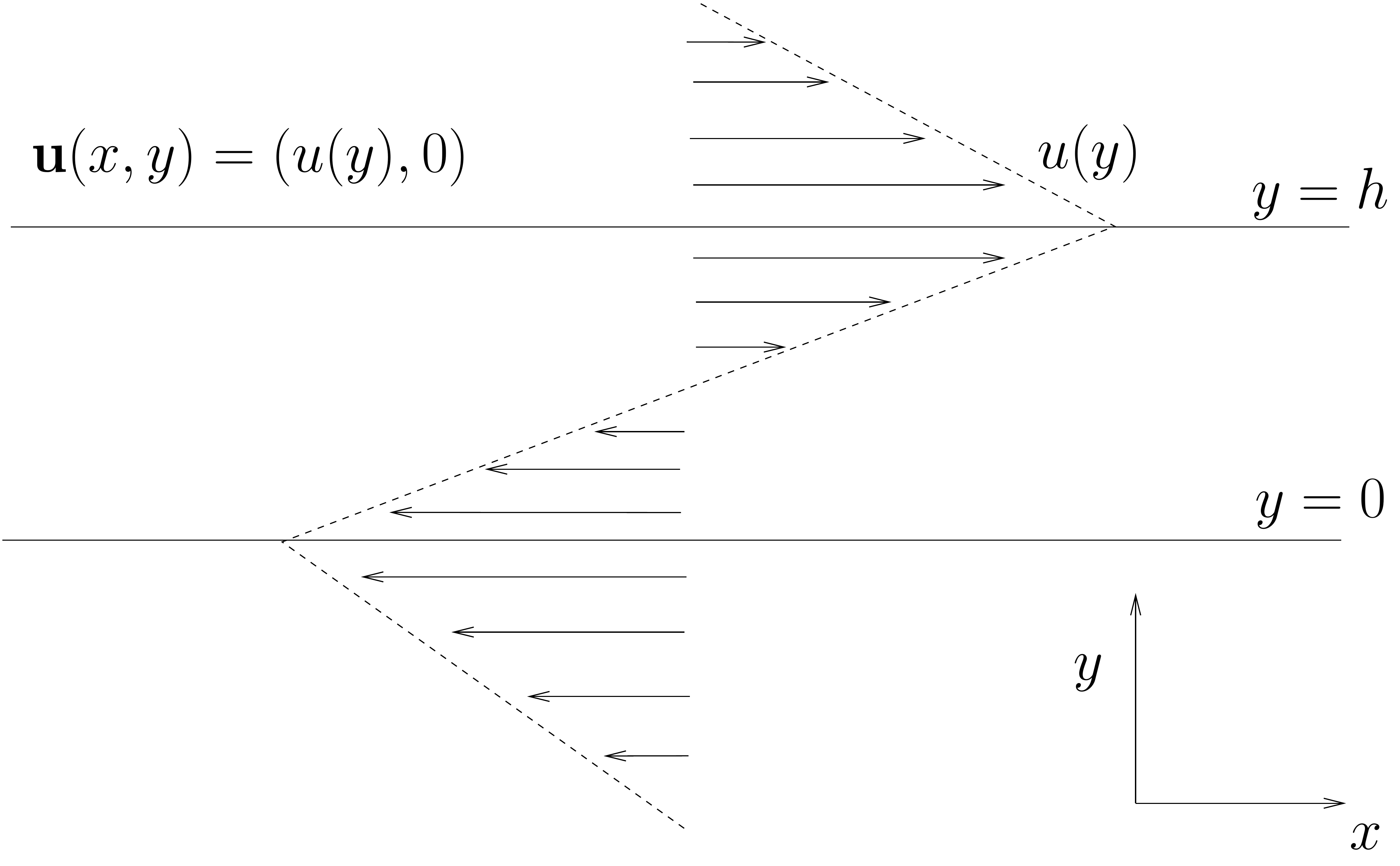}}
         \caption{Test problem of two parallel lines pulled in opposite directions by tether points.}
         \label{fig:test1_schem2}
         \end{center}
\end{figure}
See Figure \ref{fig:test1_schem2} for an illustration. We implement this scenario by using tether points, i.e.\ we prescribe the motion of virtual lines $\mb{Z}_0(t)$ and $\mb{Z}_h(t)$ and connect immersed boundaries to these lines by stiff springs with spring constant $k_\text{tether}$. Tether points are widely used in immersed boundary simulations to prescribe either rigid boundaries, as in Wiens et al.~\cite{wiens2015simulating}, or moving boundaries with prescribed kinematics, as in the simulations of jellyfish in Hamlet et al.~\cite{hamlet2011numerical}. The Lagrangian forces are given by
\begin{align}
 \mb{F}_h &= -k_\text{tether}(\mb{X}_h-\mb{Z}_h(t)) \\
  \mb{F}_0 &= -k_\text{tether}(\mb{X}_0-\mb{Z}_0(t)).
\end{align}

In our simulations we use a periodic domain of size $L_x = L_y = 2$, with a prescribed shear rate of $\dot{\gamma} = 0.15$ and a height $h = 1/24$ between the two lines. The fluid has a density of $\rho = 1$ and a dynamic viscosity of $\mu = 0.02$, which together with the characteristic velocity of $0.004$ gives a Reynolds number of $\text{Re} \approx 0.2$. The units are taken to be dimensionless. The timestep is refined according to $\Delta t = 0.01 \Delta x$. Each line is parameterized in the Lagrangian domain $[0,2\pi]$ and discretized using $N_r$ points so that $\Delta q = 2\pi/N_r$. We choose $N_r = 2N_x$ to satisfy the rule of thumb that, to prevent leakage through boundaries, the spacing between points on immersed boundaries should be approximately half the grid spacing \cite{acta}. Further, we increase the tethering constant together with the grid resolution according to $k_\text{tether} = 12.5\cdot(N_r/16)^2$, so that the lines become rigid in the limit $\Delta x \to 0$. Simulations are run up to a total time of $T = 10$.

\begin{figure}[!ht]
	\subfloat[$x$-velocity]{\label{figtest1_ib:vel}
	 \includegraphics[width=2.25 in]{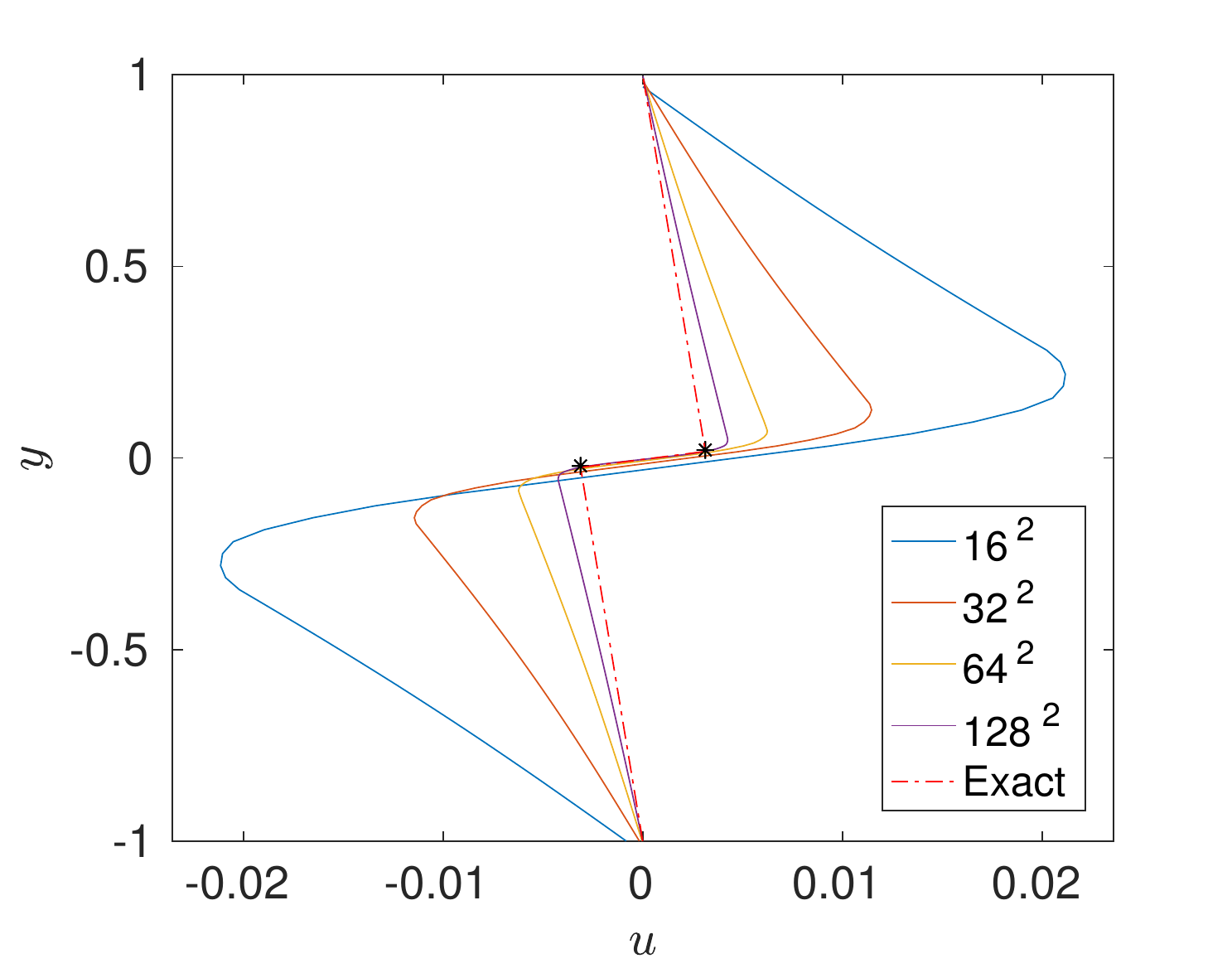}}
	\subfloat[Normal derivative of $x$-velocity]{
	 \includegraphics[width=2.25 in]{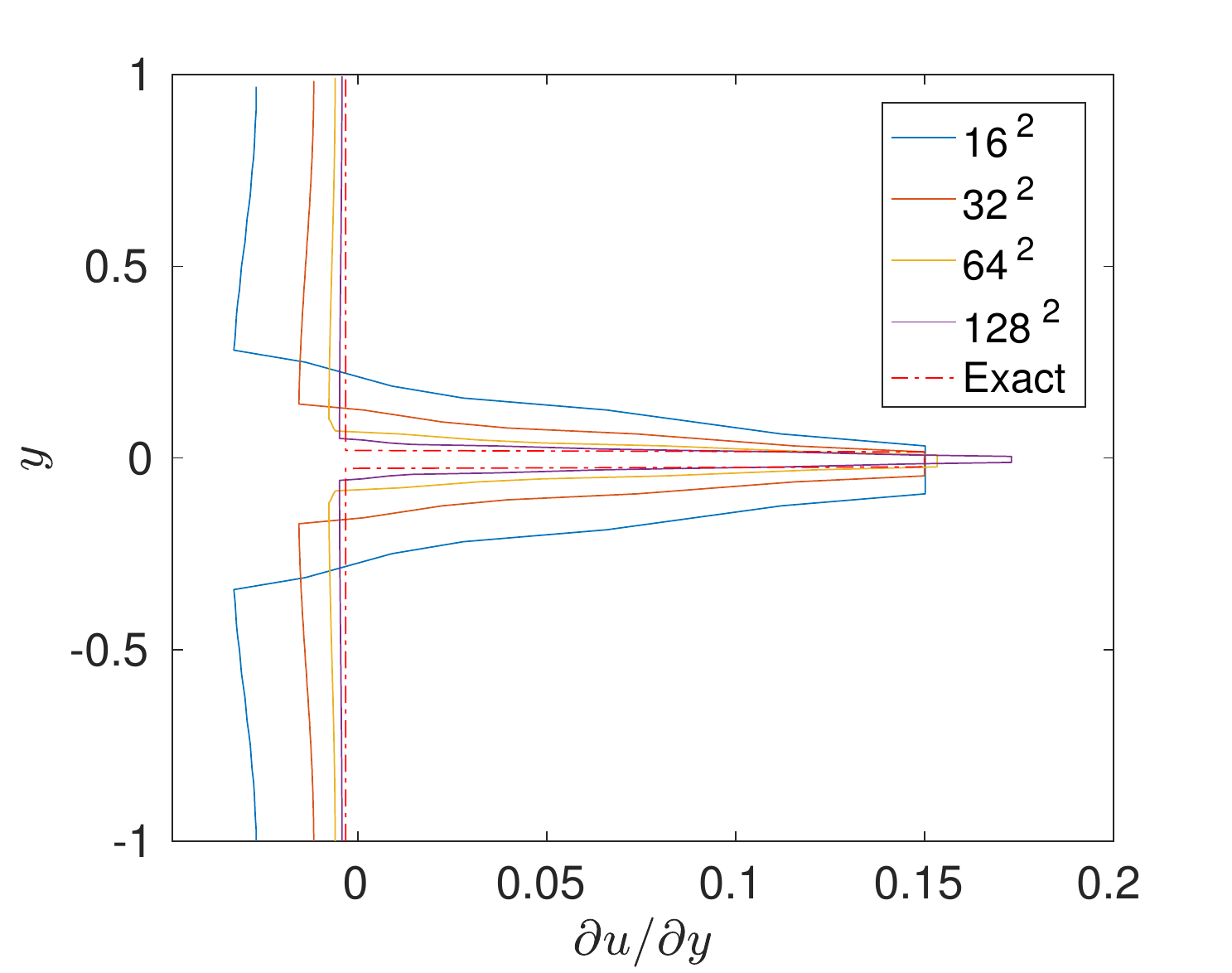}}\\
	 \subfloat[Error in $x$-velocity]{
	 \includegraphics[width=2.25 in]{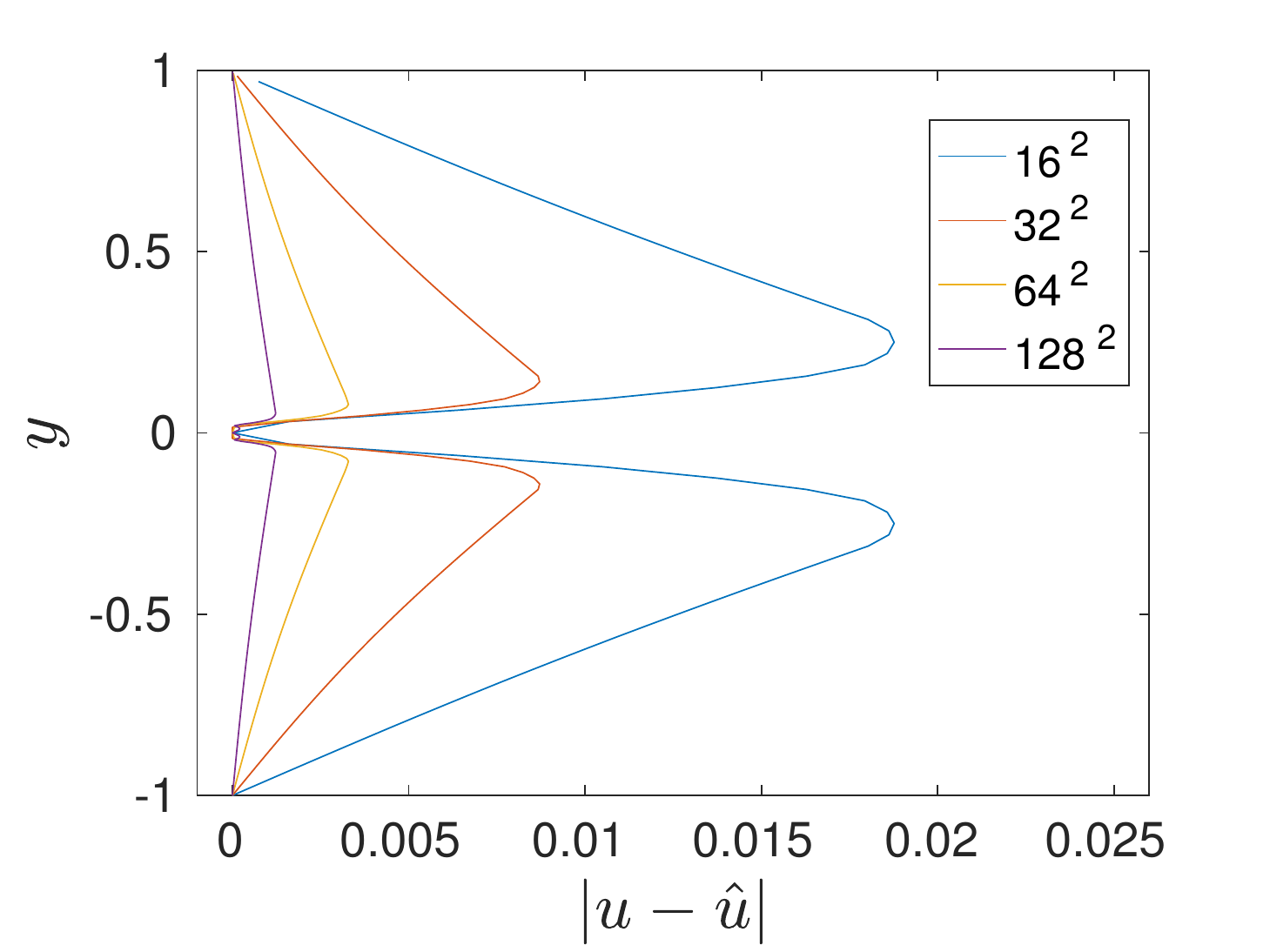}}
	\subfloat[Convergence]{\label{figtest1_ib:l1err}
	 \includegraphics[width=2.25 in]{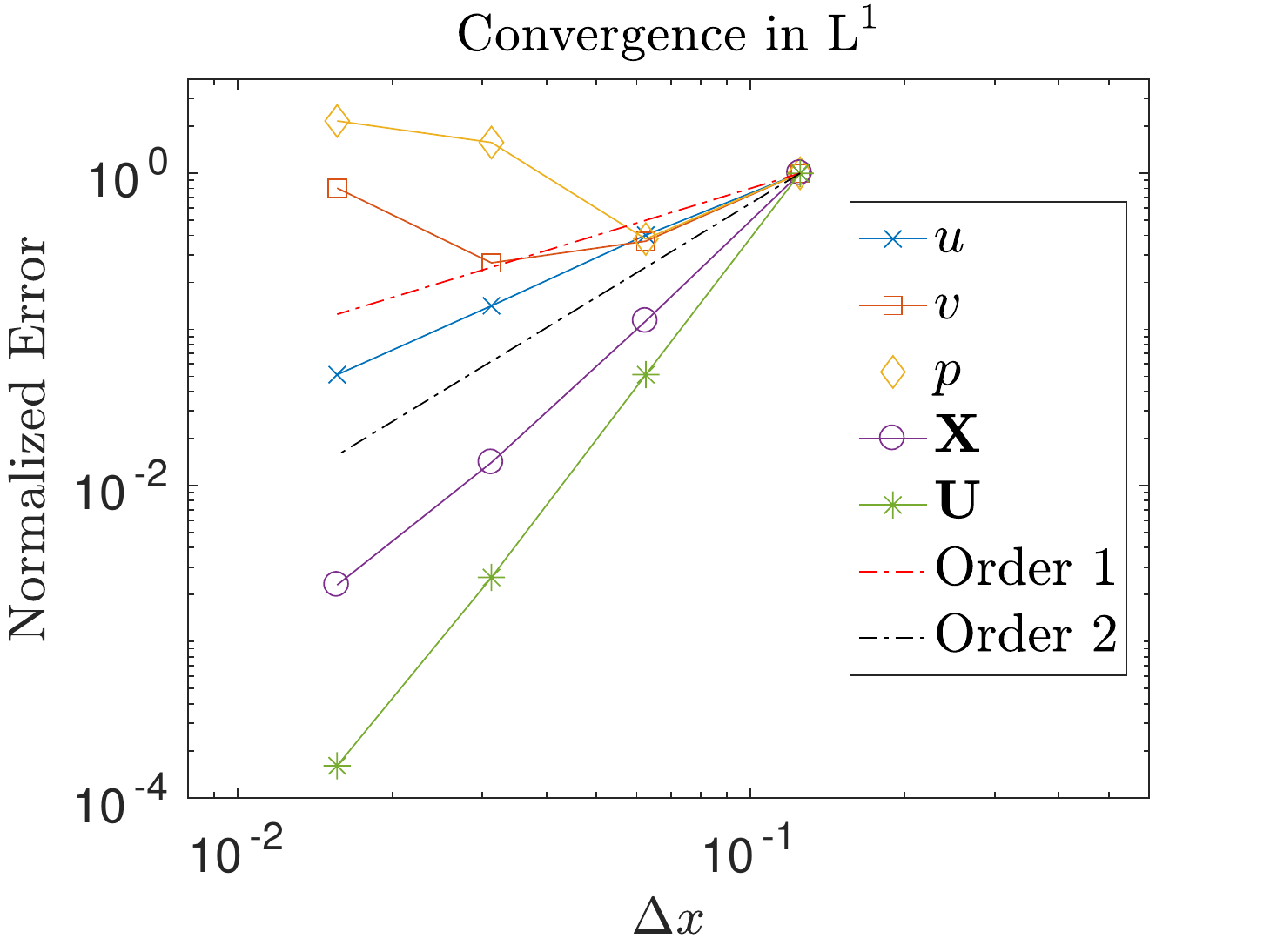}}
         \caption{Results from applying the standard immersed boundary method to the problem of parallel lines pulled in opposite directions. The asterisks in \protect\subref{figtest1_ib:vel} show the velocity values used to advect the lines. Note in \protect\subref{figtest1_ib:l1err} that all errors in $v$ and $p$ are $\lesssim 10^{-12}$, so that rounding errors dominate.}
         \label{fig:test1_ib}
\end{figure}
Given the tether point kinematics prescribed in our case, as $k_\text{tether} \to \infty$ we expect to recover $\mb{X}_0 = \mb{Z}_0$ and $\mb{X}_h = \mb{Z}_h$ with a linear shear flow throughout the domain. The results obtained using the standard immersed boundary method are illustrated in Figure \ref{fig:test1_ib}, with \edone{the mean velocities of the immersed boundaries denoted by asterisks}. The standard method yields accurate results provided that $h \gg \Delta y$. However, as would be reasonably expected it \edtwo{becomes less accurate when $h \lesssim \Delta y$}. Because of the underresolved relative velocity between the two lines, extreme spring forces are generated that lead to large overshoots in the fluid velocities near boundaries. See Table \ref{tab:test1_ib} for the errors obtained on different grid sizes, \edone{including the corresponding relative error in the mean shear rate $\dot{\gamma} = (\overline{U}_h-\overline{U}_0)/\overline{h}$, where $\overline{U}_h$ and $\overline{U}_0$ are the average velocities on each surface and $\overline{h}$ is the average height.}  

\begin{table}[]
\subfloat[$\textrm{L}^\infty$ errors]{\label{flat_ib:Linf}
\begin{tabular}{c c c c c c}
\hline
\rule{0pt}{4ex}
gridsize &$u$ &$v$ &$p$  &$\mb{X}$ & $\mb{U}$\\
\hline
\rule{0pt}{4ex}$16$  & $2.2\cdot 10^{-2}$ & $2.5\cdot 10^{-13}$ & $4.7\cdot 10^{-13}$ & $7.5\cdot 10^{-4}$ & $4.7\cdot 10^{-6}$\\
$32$  & $9.4\cdot 10^{-3}$ & $1.3\cdot 10^{-13}$ & $4.0\cdot 10^{-13}$ & $8.4\cdot 10^{-5}$ & $2.4\cdot 10^{-7}$\\
$64$  & $3.3\cdot 10^{-3}$ & $2.6\cdot 10^{-13}$ & $8.1\cdot 10^{-12}$ & $1.0\cdot 10^{-5}$ & $1.2\cdot 10^{-8}$\\
$128$  & $1.2\cdot 10^{-3}$ & $8.9\cdot 10^{-13}$ & $1.3\cdot 10^{-11}$ & $1.7\cdot 10^{-6}$ & $7.5\cdot 10^{-10}$
\end{tabular}
}\\
\subfloat[$\textrm{L}^1$ errors and relative error in $\dot{\gamma}$.]{\label{flat_ib:L1}
\begin{tabular}{c c c c c c c}
\hline
\rule{0pt}{4ex}
gridsize &$u$ &$v$ &$p$ &$\mb{X}$ & $\mb{U}$ &$\dot{\gamma}$\\
\hline
\rule{0pt}{4ex}$16$  & $4.4\cdot 10^{-2}$ & $2.2\cdot 10^{-13}$ & $1.3\cdot 10^{-13}$ & $4.7\cdot 10^{-3}$ & $3.0\cdot 10^{-5}$ & $1.5\cdot 10^{-3}$\\
$32$  & $1.8\cdot 10^{-2}$ & $8.0\cdot 10^{-14}$ & $4.9\cdot 10^{-14}$ & $5.3\cdot 10^{-4}$ & $1.5\cdot 10^{-6}$ & $7.7\cdot 10^{-5}$\\
$64$  & $6.3\cdot 10^{-3}$ & $5.8\cdot 10^{-14}$ & $2.0\cdot 10^{-13}$ & $6.6\cdot 10^{-5}$ & $7.6\cdot 10^{-8}$ & $3.9\cdot 10^{-6}$\\
$128$  & $2.3\cdot 10^{-3}$ & $1.8\cdot 10^{-13}$ &  $2.8\cdot 10^{-13}$ & $1.1\cdot 10^{-5}$ & $4.7\cdot 10^{-9}$ & $2.4\cdot 10^{-7}$\\
\end{tabular}
}
\caption{Convergence results for the standard immersed boundary method applied to the problem of parallel lines pulled in opposite directions. We give \edone{absolute} $\textrm{L}^1$ and $\textrm{L}^\infty$ errors in the velocity components $(u,v)$, pressure $p$, and immersed boundary position $\mb{X}$ and velocity $\mb{U}$. We also report the \edone{relative} error in the mean shear rate $\dot{\gamma}$.}
\label{tab:test1_ib}
\end{table}

\begin{figure}[!ht]
	\subfloat[$x$-velocity]{\label{figtest1_lube:vel}
	 \includegraphics[width=2.25 in]{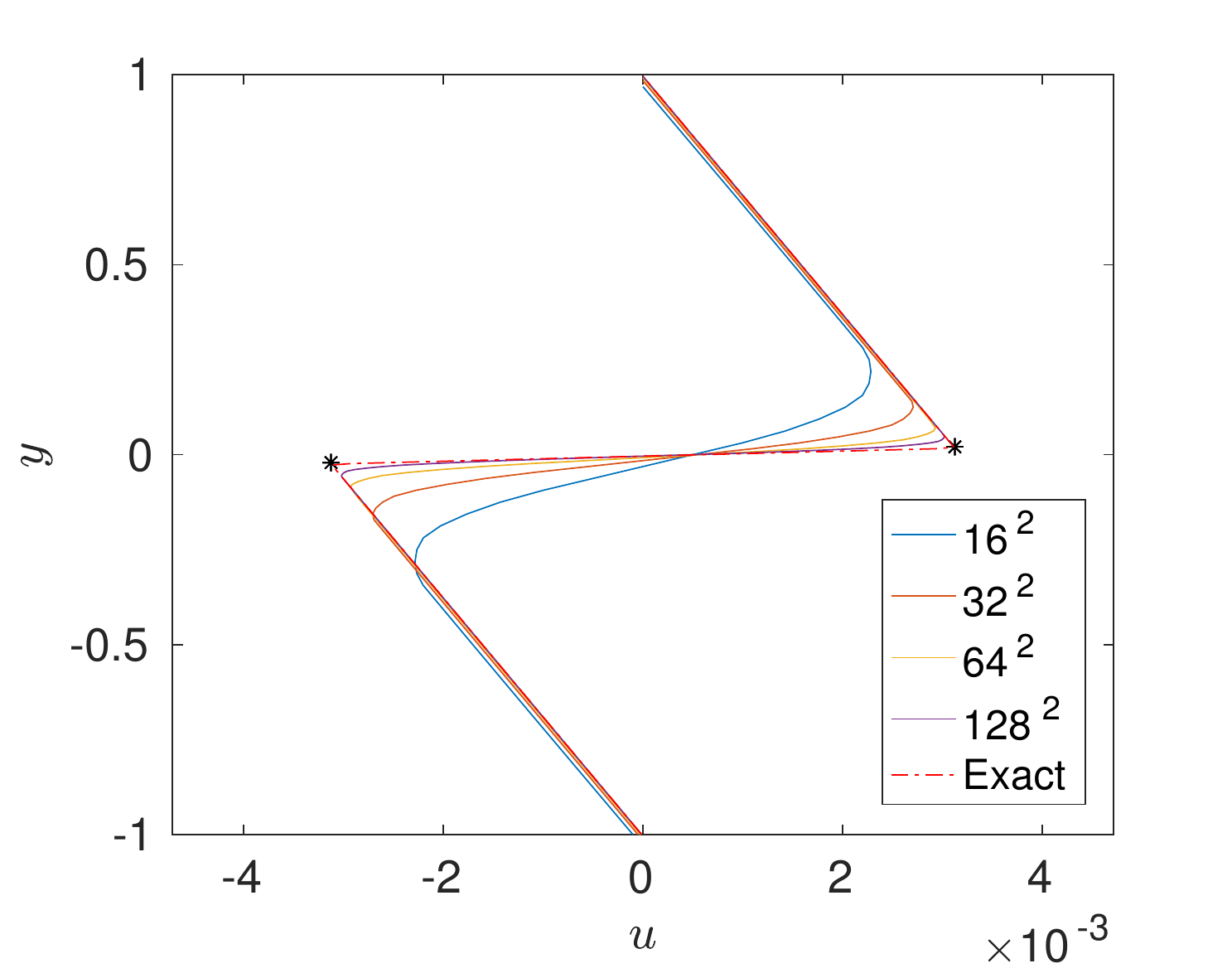}}
	\subfloat[Normal derivative of $x$-velocity]{
	 \includegraphics[width=2.25 in]{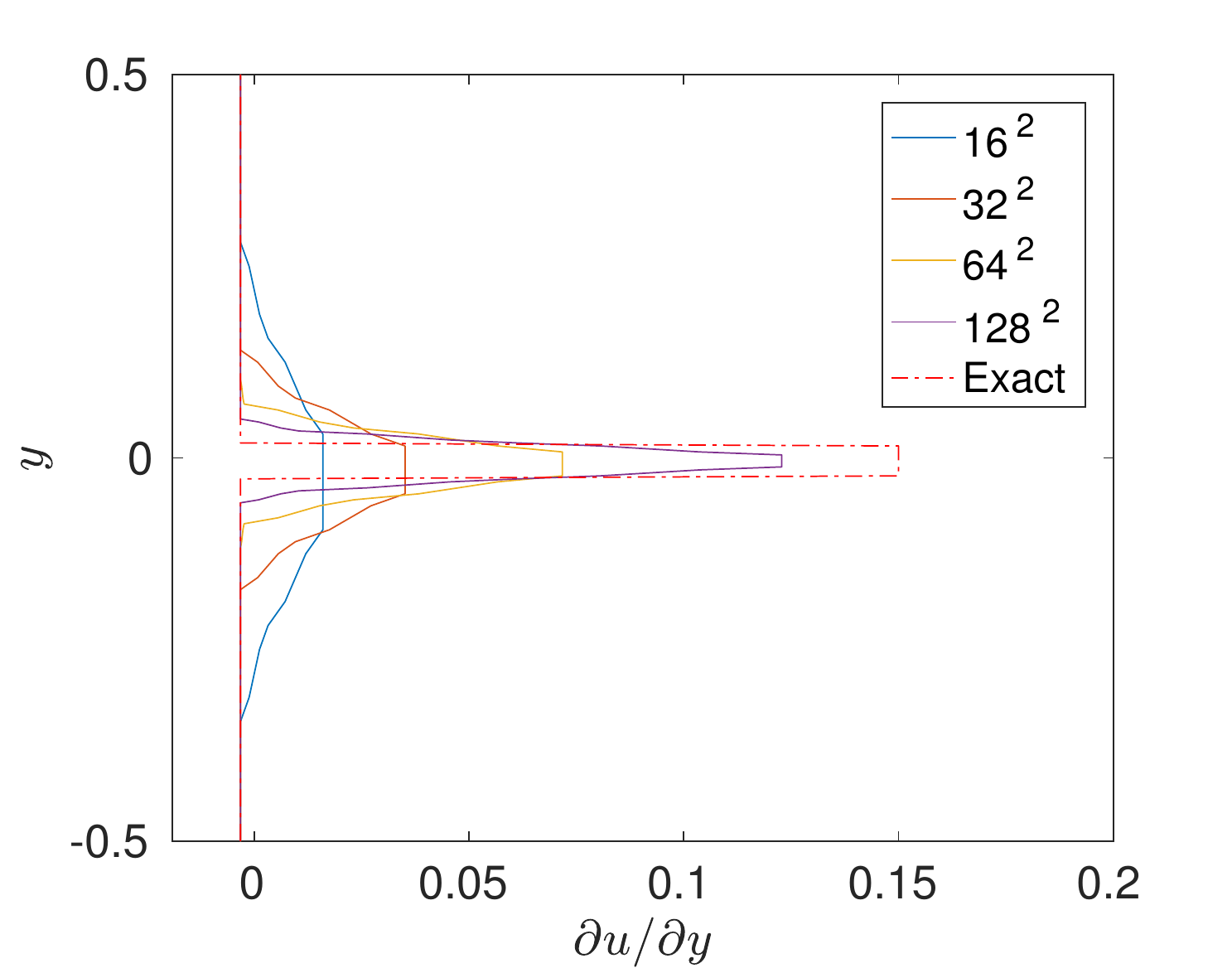}}\\
	 \subfloat[Error in $x$-velocity]{
	 \includegraphics[width=2.25 in]{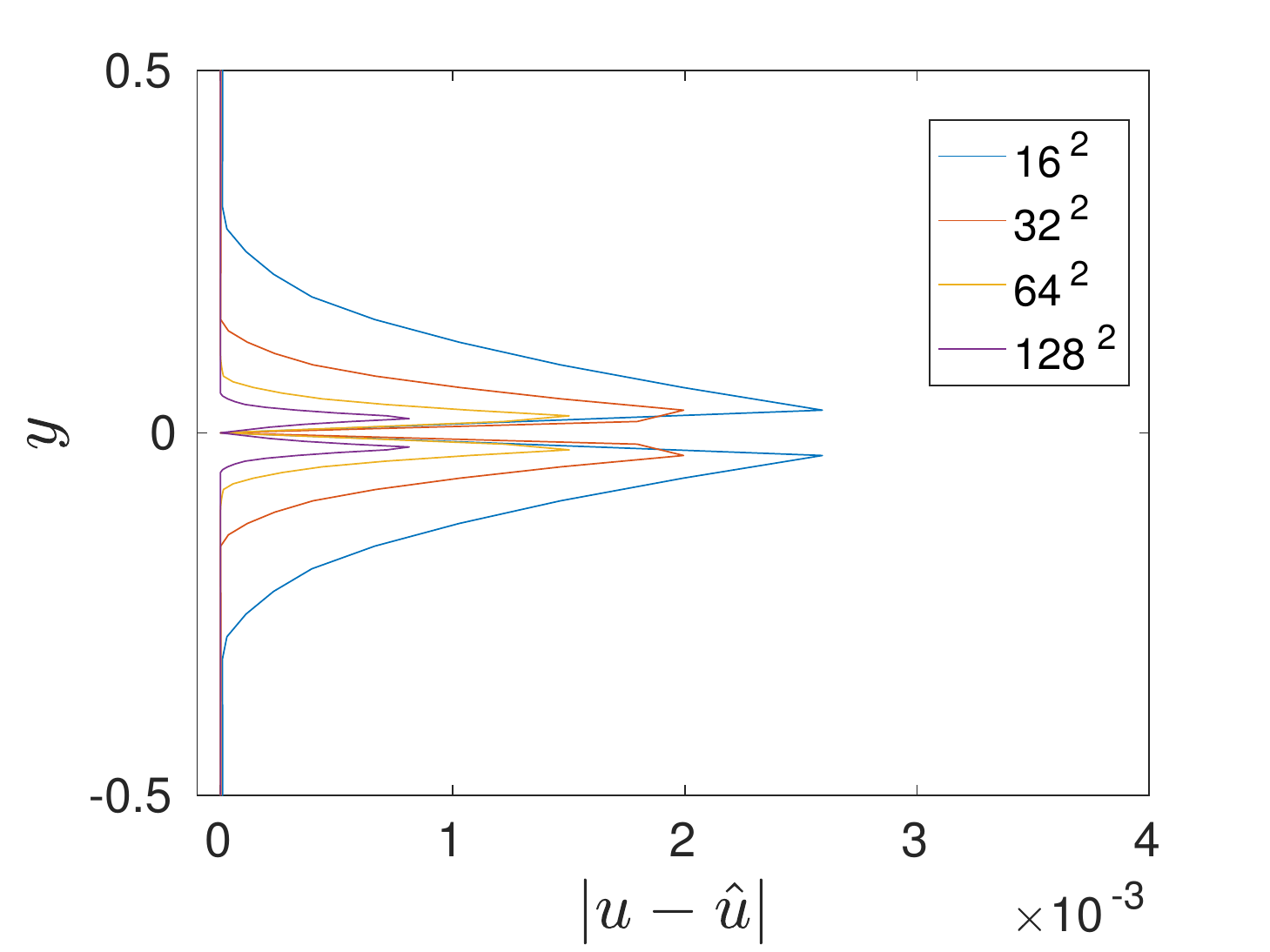}}
	\subfloat[Convergence]{\label{figtest1_lube:l1err}
	 \includegraphics[width=2.25 in]{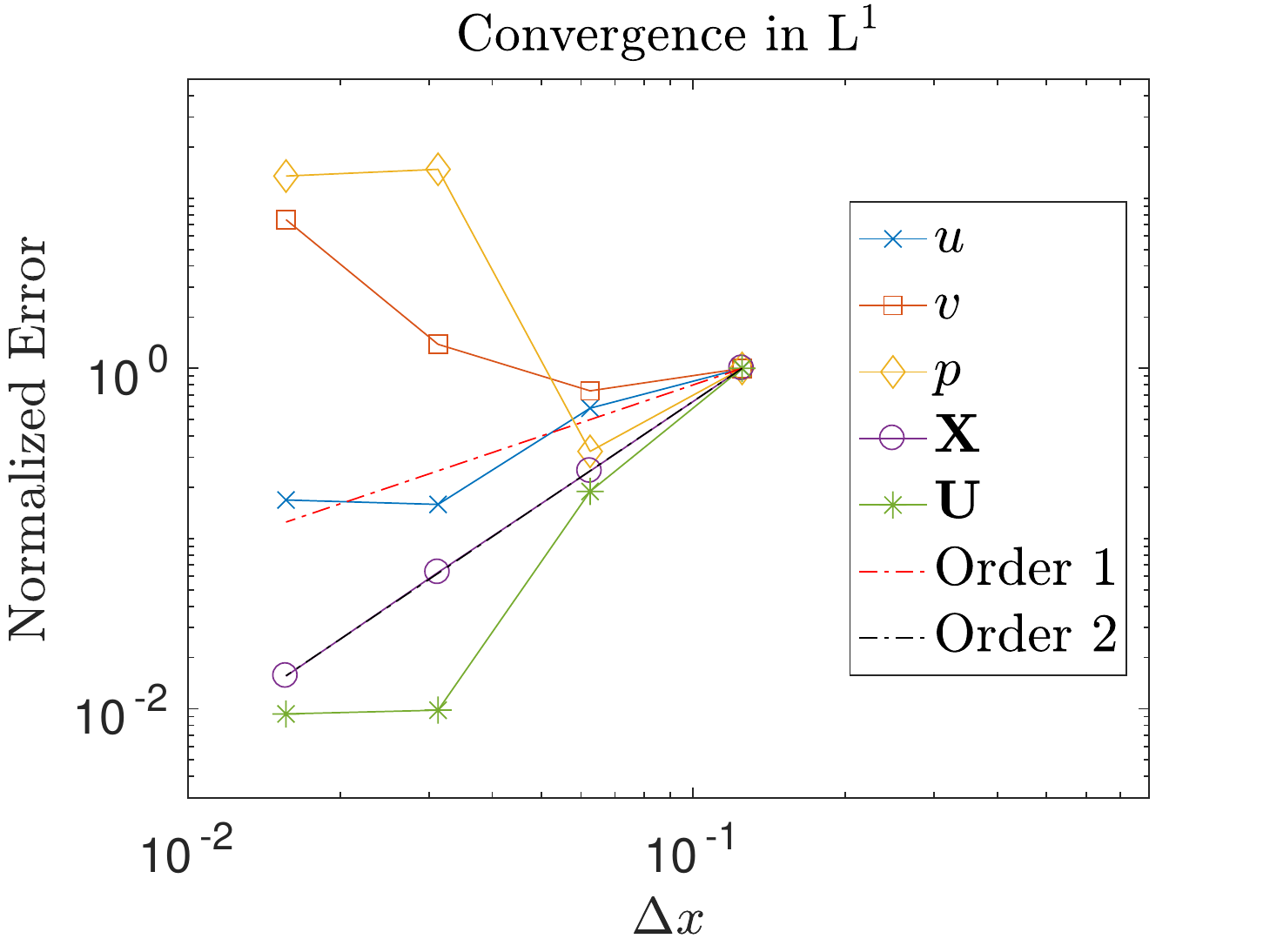}}
         \caption{Results from applying the lubricated immersed boundary method to the problem of parallel lines pulled in opposite directions. The asterisks in \protect\subref{figtest1_lube:vel} show the velocity values used to advect the lines. Because the lubricated immersed boundary method is used, the asterisks need not coincide with the interpolated velocities. Note in \protect\subref{figtest1_lube:l1err} that all errors in $v$ and $p$ are $\lesssim 10^{-13}$, so that rounding errors dominate.}
         \label{fig:test1_lube}
\end{figure}

The lubricated immersed boundary method, on the other hand, resolves the fluid flow through the thin gap exactly in the special case of linear shear and therefore converges rapidly to the correct velocities. As shown in Figure \ref{fig:test1_lube}, the resulting fluid velocities do not exhibit the large overshoots that appear in the standard immersed boundary method. The errors computed by the lubricated immersed boundary method are significantly smaller. For instance, on a coarse $16 \times 16$ grid for which $\Delta y = 3h$, the errors made by the lubricated method are over an order of magnitude smaller (see Table \ref{tab:test1_lube} for $\textrm{L}^1$ and $\textrm{L}^\infty$ errors on different grid sizes and Movies \texttt{flat\_ib64.gif} and \texttt{flat\_lube64.gif} for a visual comparison of the two methods). \edone{Because lubrication corrections are applied in this case, the asterisks denoting the mean velocities of the immersed boundaries do not coincide with the interpolated velocities from the fluid grid.}

On coarse grids that underresolve the fluid, using the lubricated immersed boundary method results in accuracies that would require a $2-4$ times finer grid if the standard immersed boundary method were used. Since the number of degrees of freedom scales with the inverse grid spacing squared, this implies that the lubricated method can obtain the same accuracy approximately $8-64$ times faster, taking  into account both the number of degrees of freedom and the CFL stability restriction $\Delta t \lesssim \Delta x$.

Note that on a $128 \times 128$ grid for which $\Delta y = 3h/8$, the standard immersed boundary method gives nearly the same error as the lubricated method. This implies that when the gap contains more than 2 gridpoints, the linear shear flow is sufficiently resolved and the lubricated method is no longer necessary.

\begin{table}[]
\subfloat[$\textrm{L}^\infty$ errors]{\label{flat_lube:Linf}
\begin{tabular}{c c c c c c}
\hline
\rule{0pt}{4ex}
gridsize &$u$ &$v$ &$p$  &$\mb{X}$ & $\mb{U}$\\
\hline
\rule{0pt}{4ex}$16$  & $1.6\cdot 10^{-3}$ & $2.5\cdot 10^{-14}$ & $4.7\cdot 10^{-14}$ & $7.7\cdot 10^{-5}$ & $3.0\cdot 10^{-9}$\\
$32$  & $1.6\cdot 10^{-3}$ & $6.9\cdot 10^{-15}$ & $2.6\cdot 10^{-14}$ & $1.9\cdot 10^{-5}$ & $5.7\cdot 10^{-10}$\\
$64$  & $8.6\cdot 10^{-4}$ & $1.2\cdot 10^{-13}$ & $1.0\cdot 10^{-12}$ & $4.9\cdot 10^{-6}$ & $2.9\cdot 10^{-11}$\\
$128$  & $4.8\cdot 10^{-4}$ & $5.0\cdot 10^{-13}$ & $3.5\cdot 10^{-12}$ & $1.2\cdot 10^{-6}$ & $2.8\cdot 10^{-11}$
\end{tabular}
}\\
\subfloat[$\textrm{L}^1$ errors and relative error in $\dot{\gamma}$.]{\label{flat_lube:L1}
\begin{tabular}{c c c c c c c}
\hline
\rule{0pt}{4ex}
gridsize &$u$ &$v$ &$p$ &$\mb{X}$ & $\mb{U}$ &$\dot{\gamma}$\\
\hline
\rule{0pt}{4ex}$16$  & $8.5\cdot 10^{-4}$  & $1.4\cdot 10^{-14}$ & $1.4\cdot 10^{-14}$ & $4.8\cdot 10^{-4}$ & $1.9\cdot 10^{-8}$ & $9.6\cdot 10^{-7}$\\
$32$  & $5.0\cdot 10^{-4}$  & $1.0\cdot 10^{-14}$ & $4.5\cdot 10^{-15}$ & $1.2\cdot 10^{-4}$ & $3.6\cdot 10^{-9}$ & $1.8\cdot 10^{-7}$\\
$64$  & $1.3\cdot 10^{-4}$  & $1.9\cdot 10^{-14}$ & $2.0\cdot 10^{-13}$ & $3.1\cdot 10^{-5}$ & $1.8\cdot 10^{-10}$ & $9.4\cdot 10^{-9}$\\
$128$  & $1.4\cdot 10^{-4}$  & $1.0\cdot 10^{-13}$ &  $1.9\cdot 10^{-13}$ & $7.5\cdot 10^{-6}$ & $1.8\cdot 10^{-10}$ & $8.9\cdot 10^{-9}$\\
\end{tabular}
}
\caption{Convergence results for the lubricated immersed boundary method applied to the problem of parallel lines pulled in opposite directions.}
\label{tab:test1_lube}
\end{table}

\subsection{Eccentric rotating cylinders}
The test problem of two parallel lines is a relatively simple case, of course. The assumptions of lubrication theory are actually valid for linear shear flow independent of the gap size, and the absence of curvature hides approximations made by the lubricated immersed boundary method. Therefore, we next consider the more demanding problem of eccentric rotating cylinders.

This problem is appealing because it is simple to describe and implement but gives rise to non-trivial fluid dynamics. Unlike the case of linear shear flow, the lubrication approximation is not exact in this case. Moreover, interesting phenomena emerge such as lift on the inner cylinder and counter-rotation of the fluid on the wide side of the gap, as we describe next. In discussing the setup of this well-known problem we closely follow the exposition of Acheson \cite{Acheson1990}.

On a domain of size $L_x=L_y=2$, we use inner and outer cylinders of radii $r_1 = 3/4$ and $r_2 = 3/4(1+1/24)$ respectively, so that the nondimensional thickness $\epsilon:=(r_2-r_1)/r_1=1/24$. The outer cylinder is centered at the origin, whereas the inner cylinder is  centered at $(x_0,0)$ with $x_0=3/128$ so that the cylinders are not concentric. With this choice of parameters, the dimensionless eccentricity $\lambda:=x_0/(r_1 \epsilon)$ satisfies $\lambda= 3/4$. This value is chosen partly because whenever $\lambda > (\sqrt{13}-3)/2 \approx 0.30$, the exact solution exhibits a counter-rotating vortex on the wide side of the gap \cite{Acheson1990}. Moreover, there is a lift $F$ on the cylinder equal to $F=12 \pi \mu U \lambda/(\epsilon^2 \sqrt{1-\lambda}(2+\lambda^2))$.

Next, we investigate which of these features can be recovered on coarse grids using the lubricated immersed boundary method. To simulate this problem, we again use tether points with specified kinematics. The outer cylinder is taken to be stationary, while the inner cylinder is tethered to points rotating counterclockwise with prescribed angular velocity $U=8.33 \cdot 10^{-4}$. The simulation is run up to a total time of $T = 100$. We use a fourfold stronger spring constant of $k_\text{tether} = 50\cdot(N_r/16)^2$ compared to the linear shear simulations of Section \ref{sec:simp_shear}. All other parameters are unchanged.

\begin{figure}[!ht]
	\subfloat[$v$ through a horizontal slice]{
	 \includegraphics[width=2.5 in]{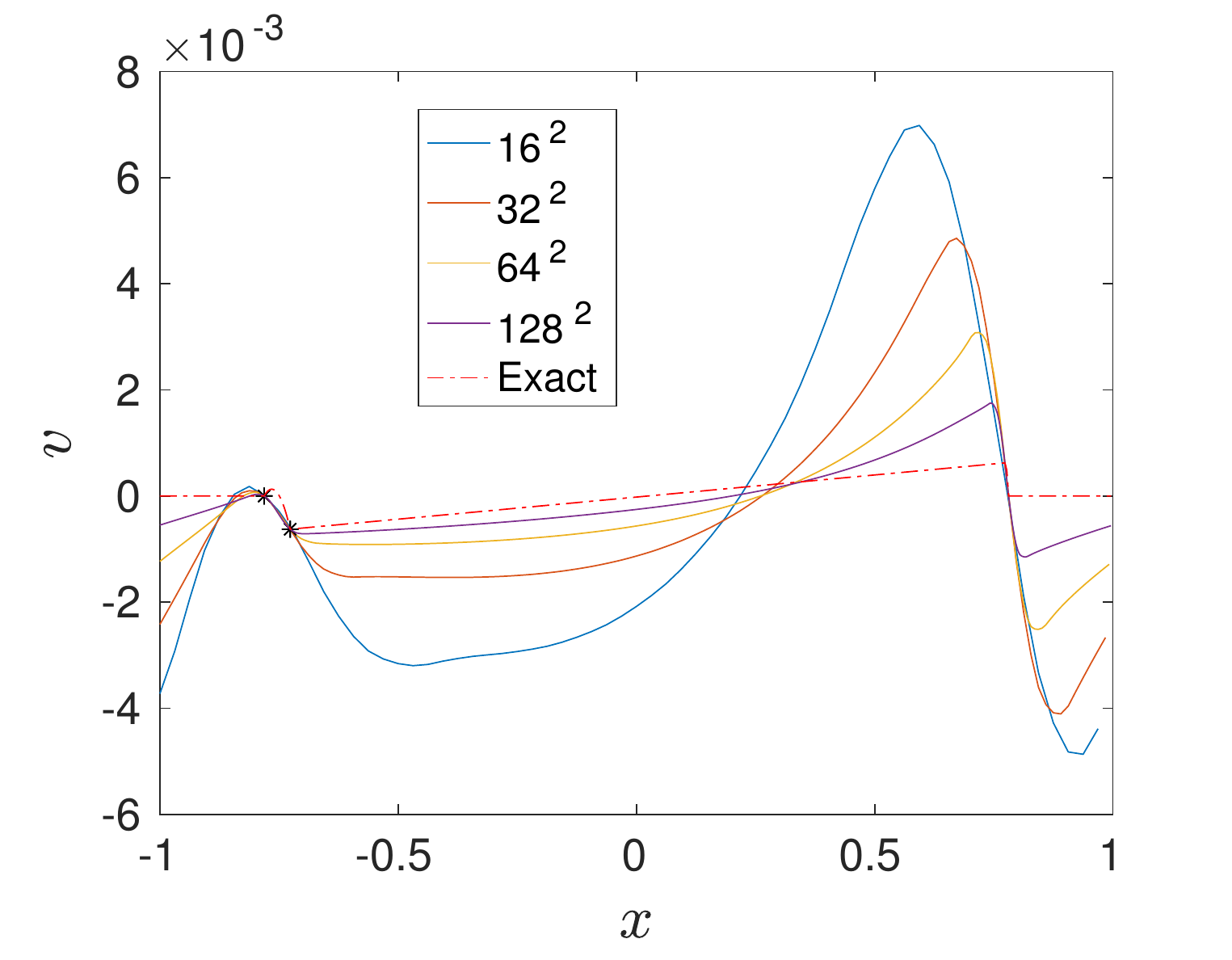}}
	\subfloat[Pressure gradient]{
	 \includegraphics[width=2.5 in]{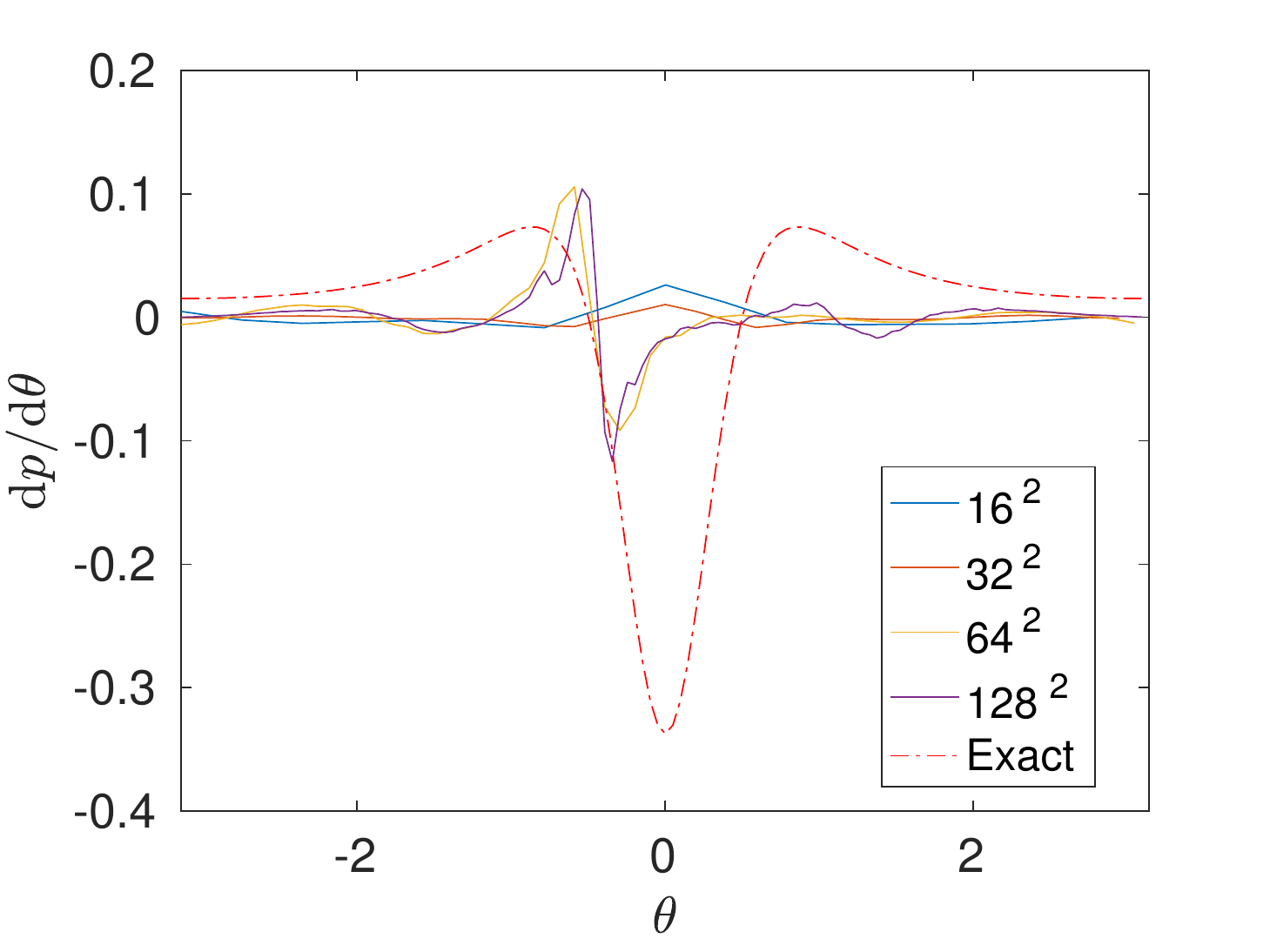}}\\
	\subfloat[Gap velocity $u_\theta$ in rescaled coordinates on a $32^2$ grid]{\label{figtest2_ib:ugap}
	 \includegraphics[width=2.5 in]{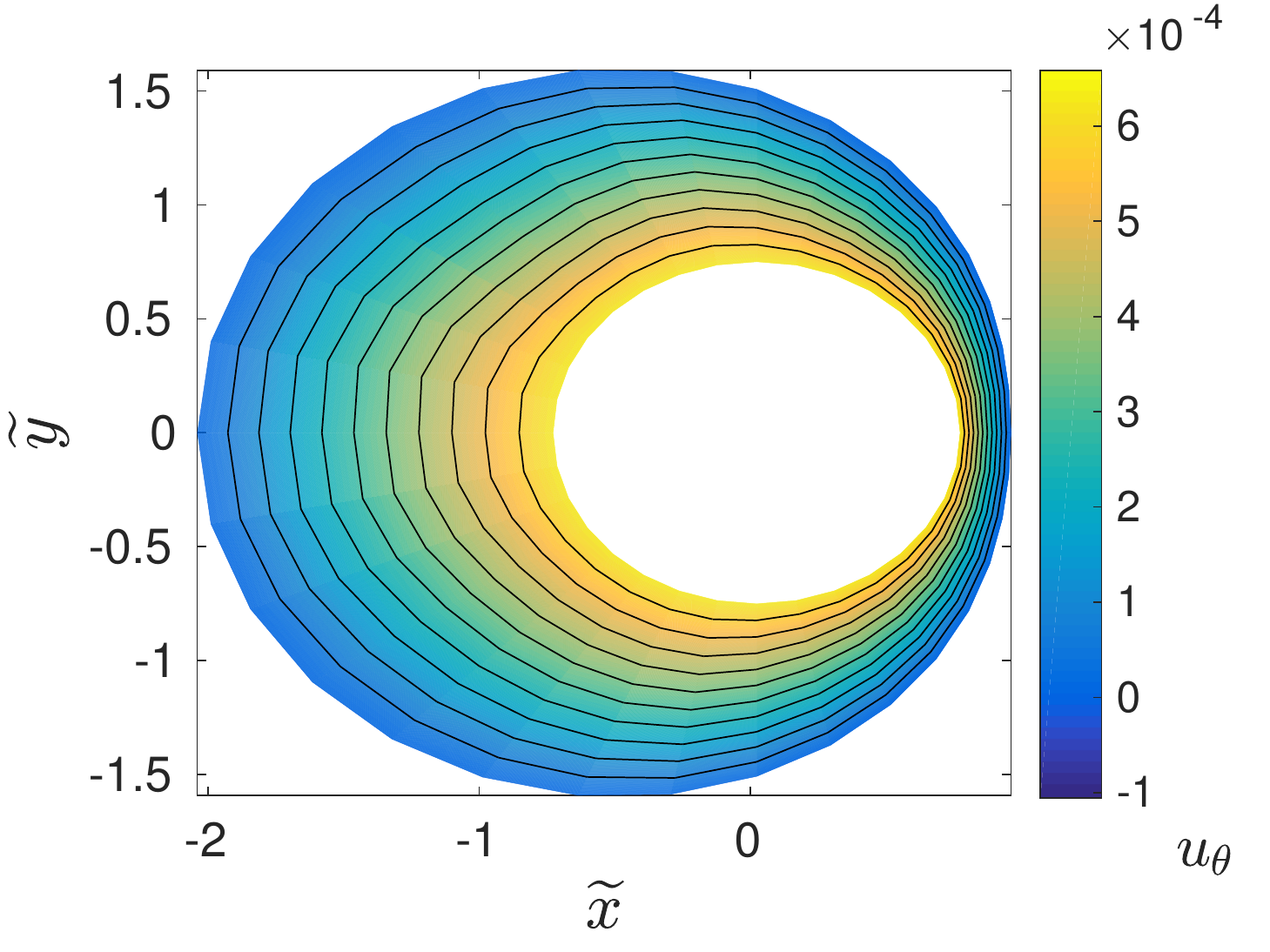}}
	 \subfloat[Convergence]{
	 \includegraphics[width=2.5 in]{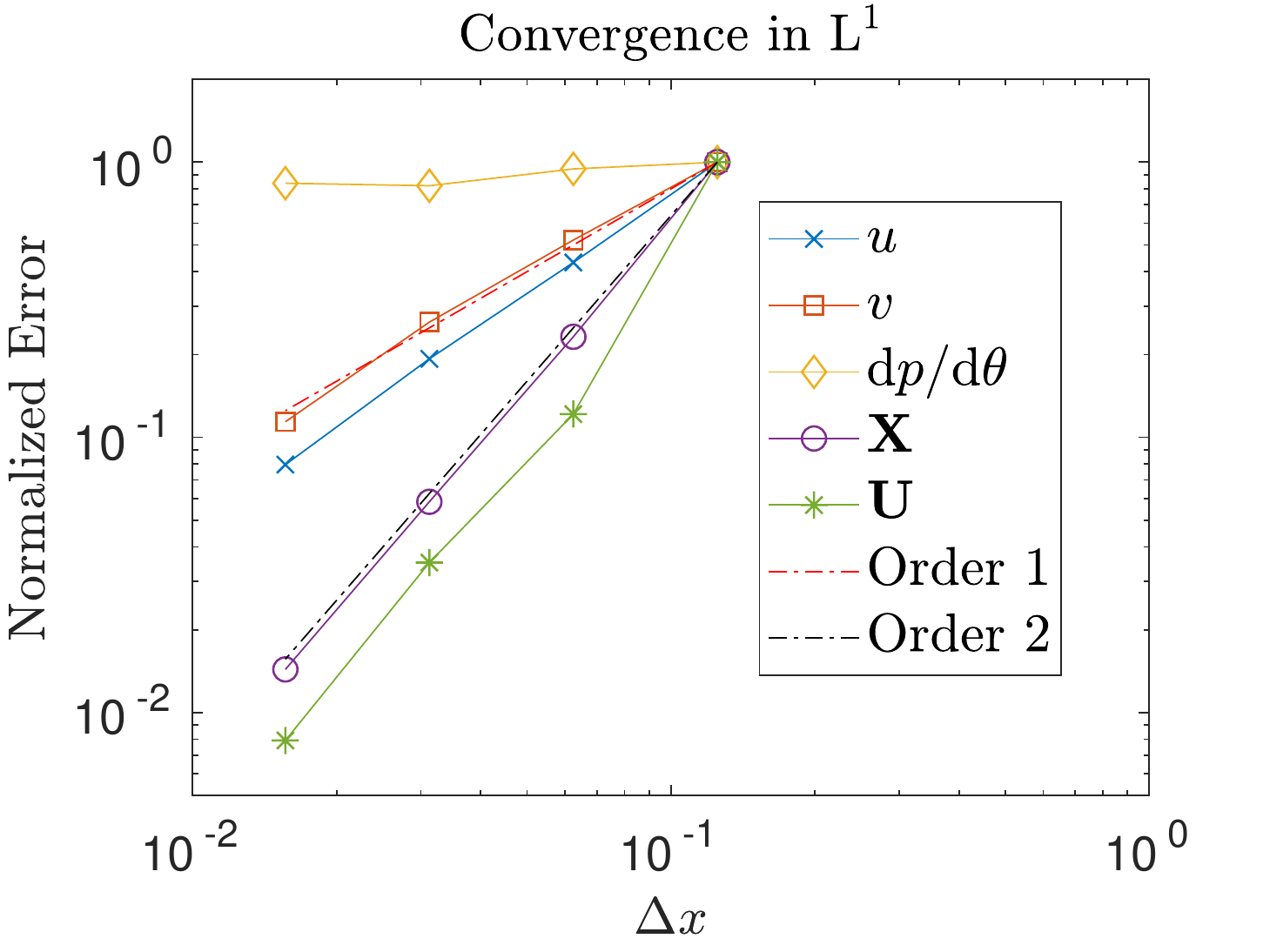}}
         \caption{Standard immersed boundary method applied to the problem of eccentric rotating cylinders. Note the absence of reverse flow in \protect\subref{figtest2_ib:ugap}, i.e.\ $u_\theta\ge0$ everywhere. The gap in \protect\subref{figtest2_ib:ugap} has been blown up for purposes of visualization by scaling all distances from the inner cylinder by a factor of $1/\epsilon = 24$.}
         \label{fig:test2_ib}
\end{figure}
\begin{figure}[!ht]
	\subfloat[$v$ through a horizontal slice]{
	 \includegraphics[width=2.5 in]{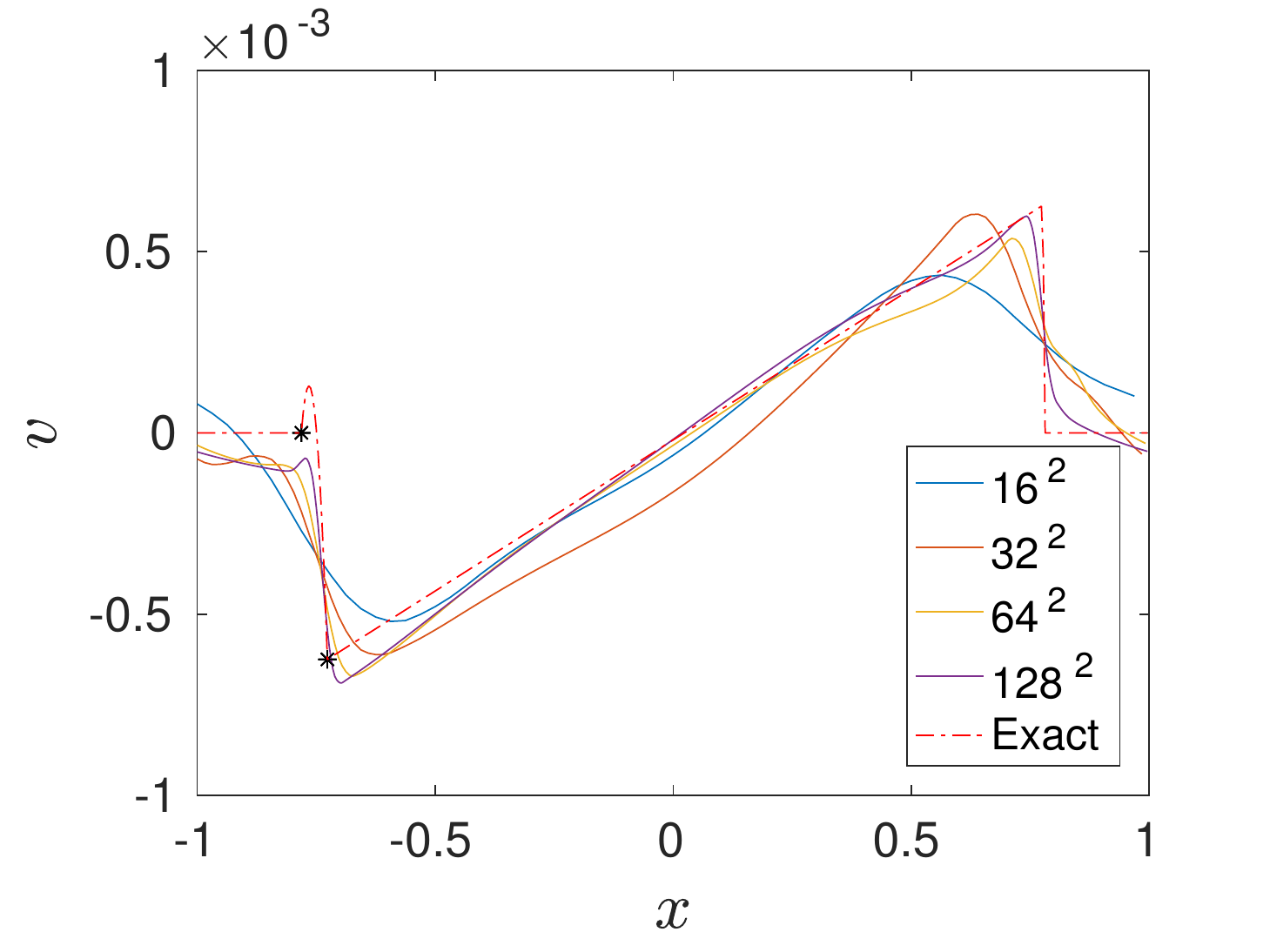}}
	\subfloat[Pressure gradient]{\label{figtest2_lube:press}
	 \includegraphics[width=2.5 in]{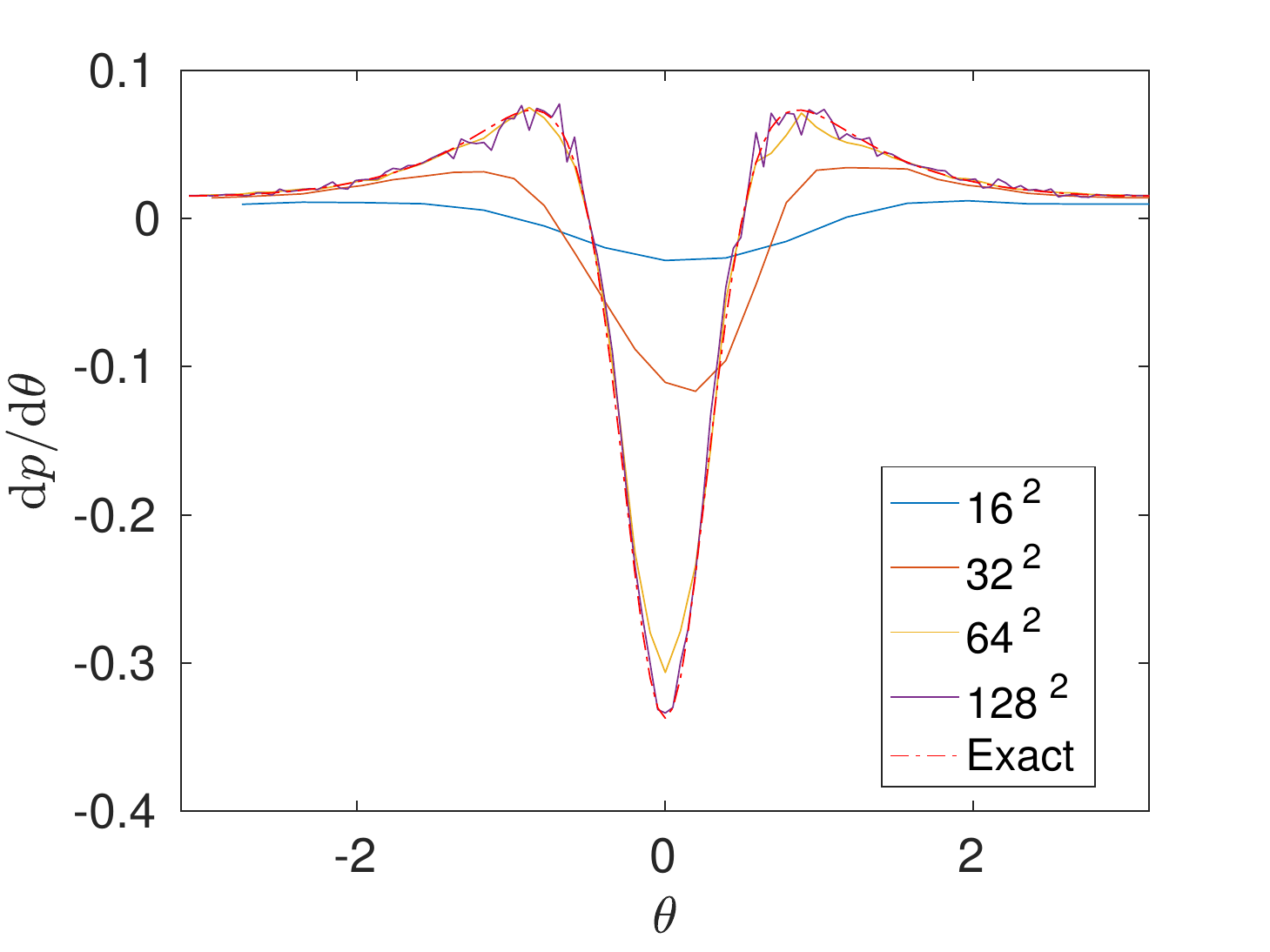}}\\
	 \subfloat[Gap velocity $u_\theta$ in rescaled coordinates on a $32^2$ grid]{\label{figtest2_lube:ugap}
	 \includegraphics[width=2.5 in]{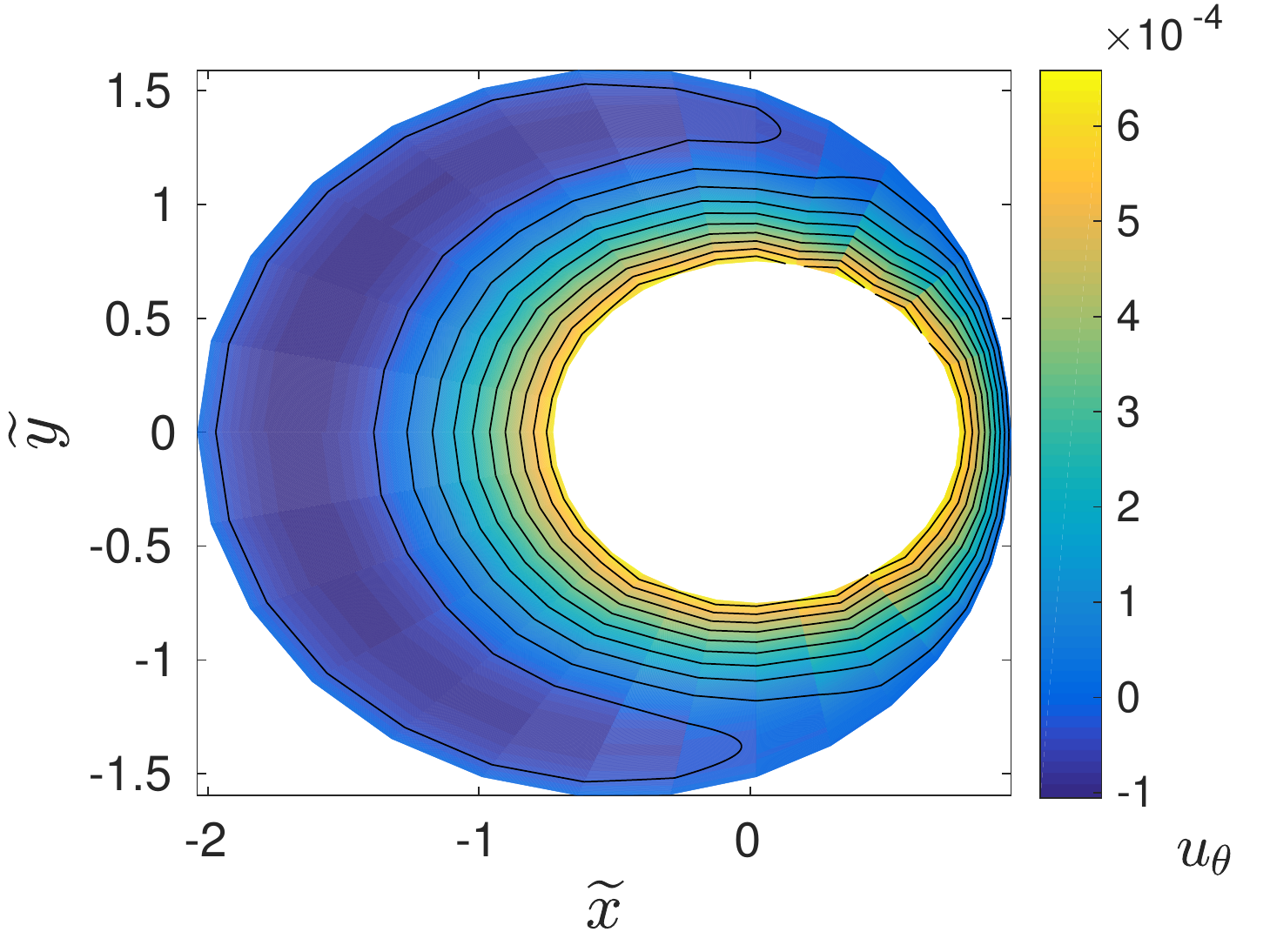}}
	 \subfloat[Convergence]{
	 \includegraphics[width=2.5 in]{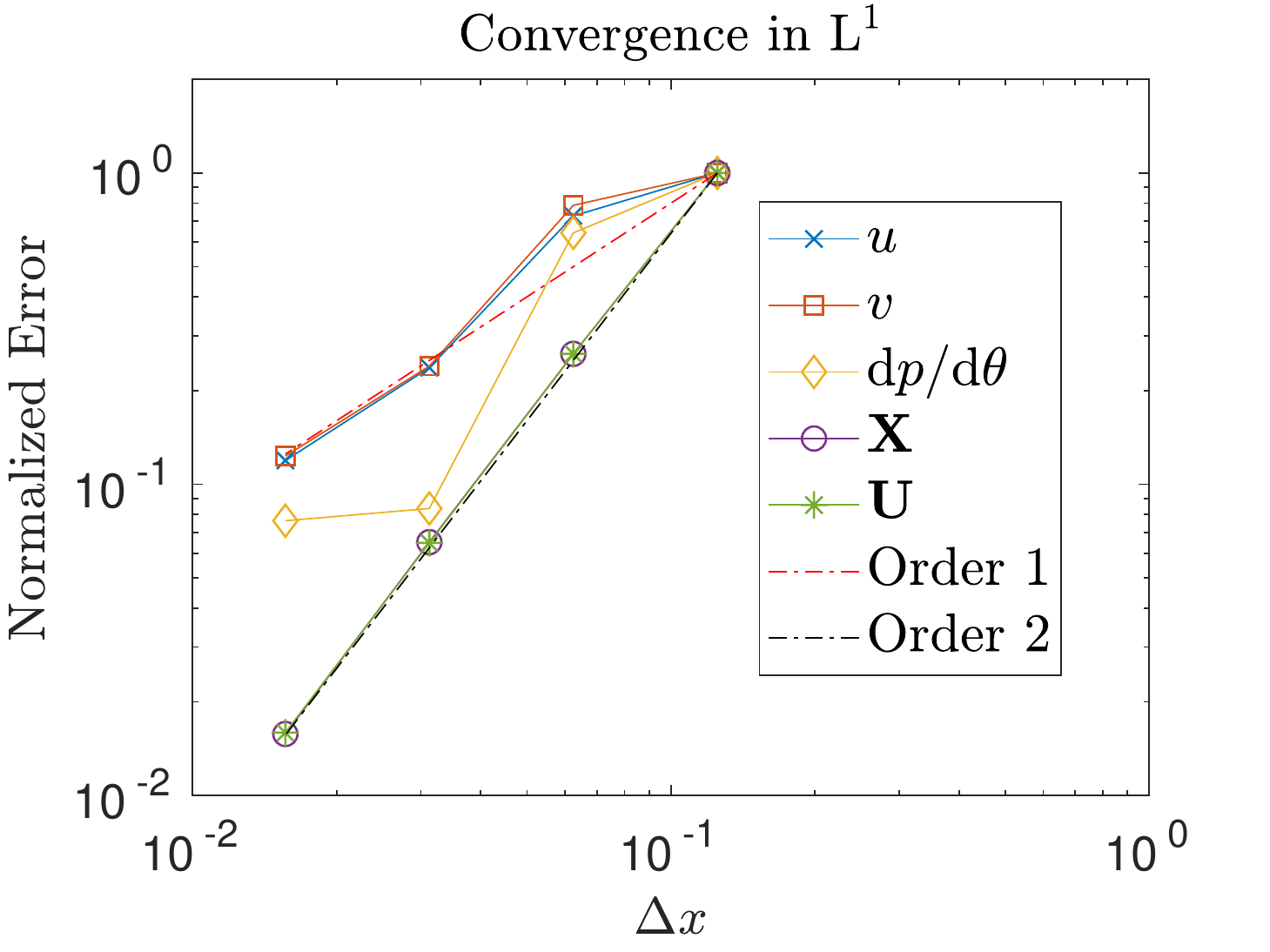}}
         \caption{Lubricated immersed boundary method applied to the problem of eccentric rotating cylinders. Note the presence of reverse flow in \protect\subref{figtest2_lube:ugap}, i.e.\ $u_\theta<0$ on the far side of the gap. The gap in \protect\subref{figtest2_lube:ugap} has been blown up for purposes of visualization by scaling all distances from the inner cylinder by a factor of $1/\epsilon = 24$.}
         \label{fig:test2_lube}
\end{figure}

On the coarsest $16 \times 16$ grid used, the two cylinders come to within a minimal distance of $\Delta x/16$, and attain a maximal distance of about $\Delta x/2$. Given these subgrid separations, it is no great surprise based on the results of the previous section that using the standard immersed boundary method leads to large overshoots near the boundaries that are significantly improved by the lubricated method (Figures \ref{fig:test2_ib} and \ref{fig:test2_lube}). Further, by computing the tangential velocity through the gap using \eqref{eq:u_lube3}, we are able on a coarse, $16 \times 16$ grid to detect reverse flow indicating the existence of a counter-rotating vortex, whereas using the standard method we cannot observe this even using a significantly finer, $128 \times 128$ grid. See Figures \ref{fig:test2_ib}\subref{figtest2_ib:ugap} and \ref{fig:test2_lube}\subref{figtest2_lube:ugap} for a comparison of the two methods on a $32 \times 32$ grid.

The lubricated method accurately resolves the pressure through the gap (Figure \ref{fig:test2_lube}\subref{figtest2_lube:press}), allowing for the accurate calculation of lift. The lift is calculated by summing the Lagrangian force density over the inner immersed boundary and multiplying by $\Delta q$ to obtain a force. The lift values obtained using the lubricated immersed boundary method on $16^2$, $32^2$, $64^2$, and $128^2$ grids have corresponding relative errors 67\%, 27\%, 2.1\%, and 0.58\%, consistent with 2\textsuperscript{nd}-order accurate convergence. The lifts calculated using the standard immersed boundary method, on the other hand, actually have the wrong sign for grids smaller than $128^2$, and converge to the correct value so slowly that there is still a $14\%$ error on a $512^2$ grid. These significant errors are not entirely surprising given that the standard method does not accurately capture discontinuities at the fluid-structure interface \cite{bg}.

Interestingly, the advantages of the lubricated method are not made readily apparent by the errors in certain other quantities calculated at several grid sizes (Tables \ref{tab:test2_ib} and \ref{tab:test2_lube}). The standard immersed boundary method is able to separate the cylinders nearly as well as the lubricated method and does not seem to exhibit a significantly underresolved relative velocity between the cylinders. Indeed, the mean shear rate is the only quantity displayed in Tables \ref{tab:test2_ib} and \ref{tab:test2_lube} that results in an error at least an order of magnitude smaller for the lubricated method.  However, as in the linear shear test problem, there are large overshoots in the fluid velocities (see Movies \texttt{ecc\_ib64.gif} and \texttt{ecc\_lube64.gif}) and other quantities such as the gap pressure are not computed accurately. This shows the importance for assessing accuracy of investigating auxiliary quantities, such as the lift and flow reversal, that are not controlled directly through the tether points. 

\begin{table}[]
\subfloat[$\textrm{L}^\infty$ errors]{\label{ecc_ib:Linf}
\begin{tabular}{c c c c c c}
\hline
\rule{0pt}{4ex}
gridsize &$u$ &$v$ &$\ud p/\ud \theta$  &$\mb{X}$ & $\mb{U}$\\
\hline
\rule{0pt}{4ex}$16$ &$5.1\cdot 10^{-3}$ & $8.2\cdot 10^{-3}$ & $3.6\cdot 10^{-1}$ & $2.0\cdot 10^{-3}$ & $1.2\cdot 10^{-5}$\\
$32$ & $2.4\cdot 10^{-3}$ & $5.2\cdot 10^{-3}$ & $3.5\cdot 10^{-1}$ & $4.1\cdot 10^{-4}$ & $5.5\cdot 10^{-7}$\\
$64$  & $1.2\cdot 10^{-3}$ & $3.0\cdot 10^{-3}$ & $3.2\cdot 10^{-1}$ & $1.2\cdot 10^{-4}$ & $3.9\cdot 10^{-7}$\\
$128$  & $5.5\cdot 10^{-4}$ & $1.2\cdot 10^{-3}$ & $3.2\cdot 10^{-1}$ & $3.0\cdot 10^{-5}$ & $8.5\cdot 10^{-8}$
\end{tabular}
}
\\
\subfloat[$\textrm{L}^1$ errors]{\label{ecc_ib:L1}
\begin{tabular}{c c c c c c c}
\hline
\rule{0pt}{4ex}
gridsize &$u$ &$v$ &$\ud p/\ud \theta$ &$\mb{X}$ & $\mb{U}$ &$\dot{\gamma}$\\
\hline
\rule{0pt}{4ex}$16$  &$6.5\cdot 10^{-3}$ & $7.3\cdot 10^{-3}$ &$6.5\cdot 10^{-2}$ & $1.0\cdot 10^{-2}$ & $1.6\cdot 10^{-5}$ & $4.7\cdot 10^{-3}$\\
$32$  & $2.8\cdot 10^{-3}$ & $3.8\cdot 10^{-3}$ & $6.1\cdot 10^{-2}$ & $2.4\cdot 10^{-3}$ & $1.9\cdot 10^{-6}$ & $2.4\cdot 10^{-5}$\\
$64$  & $1.3\cdot 10^{-3}$ & $1.9\cdot 10^{-3}$ & $5.3\cdot 10^{-2}$ & $6.0\cdot 10^{-4}$ & $5.6\cdot 10^{-7}$ & $9.7\cdot 10^{-6}$\\
$128$  & $5.2\cdot 10^{-4}$ &  $8.3\cdot 10^{-4}$ & $5.5\cdot 10^{-2}$ & $1.5\cdot 10^{-4}$ & $1.3\cdot 10^{-7}$ & $1.3\cdot 10^{-6}$\\
\end{tabular}
}
\caption{Convergence results for the standard immersed boundary method applied to the problem of eccentric rotating cylinders. We give \edone{absolute} $\textrm{L}^1$ and $\textrm{L}^\infty$ errors in the velocity components $(u,v)$, pressure gradient $\ud p/\ud \theta$, and immersed boundary position $\mb{X}$ and velocity $\mb{U}$, along with the \edone{relative} error in the mean shear rate $\dot{\gamma}$ computed in terms of the angular velocities. The pressure gradient is computed by applying the jump condition \eqref{eq:jump2} to obtain $p$ and then taking a centered finite difference along the boundary.}
\label{tab:test2_ib}
\end{table}

\begin{table}[]
\subfloat[$\textrm{L}^\infty$ errors]{\label{ecc_lube:Linf}
\begin{tabular}{c c c c c c}
\hline
\rule{0pt}{4ex}
gridsize &$u$ &$v$ &$\ud p/\ud \theta$ &$\mb{X}$ & $\mb{U}$\\
\hline
\rule{0pt}{4ex}$16$ & $2.4\cdot 10^{-3}$ & $3.0\cdot 10^{-3}$ & $3.1\cdot 10^{-1}$ & $1.7\cdot 10^{-3}$ & $1.4\cdot 10^{-6}$\\
$32$  & $4.8\cdot 10^{-3}$ & $4.4\cdot 10^{-3}$ & $2.3\cdot 10^{-1}$ & $5.8\cdot 10^{-4}$ & $4.9\cdot 10^{-7}$\\
$64$ & $3.3\cdot 10^{-3}$ & $2.8\cdot 10^{-3}$ & $3.2\cdot 10^{-2}$ & $1.8\cdot 10^{-4}$ & $1.7\cdot 10^{-7}$\\
$128$  & $2.2\cdot 10^{-3}$ & $2.0\cdot 10^{-3}$ & $2.2\cdot 10^{-2}$ & $4.5\cdot 10^{-5}$ & $4.3\cdot 10^{-8}$
\end{tabular}
}
\\
\subfloat[$\textrm{L}^1$ errors]{\label{ecc_lube:L1}
\begin{tabular}{c c c c c c c}
\hline
\rule{0pt}{4ex}
gridsize &$u$ &$v$ &$\ud p/\ud \theta$ &$\mb{X}$ & $\mb{U}$ &$\dot{\gamma}$\\
\hline
\rule{0pt}{4ex}$16$  & $1.8\cdot 10^{-3}$ & $1.9\cdot 10^{-3}$ & $5.0\cdot 10^{-2}$ & $9.1\cdot 10^{-3}$ & $7.7\cdot 10^{-6}$ & $5.4\cdot 10^{-6}$\\
$32$  & $1.3\cdot 10^{-3}$ & $1.5\cdot 10^{-3}$ & $3.2\cdot 10^{-2}$ & $2.4\cdot 10^{-3}$ & $2.0\cdot 10^{-6}$ & $4.9\cdot 10^{-7}$\\
$64$  & $4.2\cdot 10^{-4}$ & $4.6\cdot 10^{-4}$ & $4.1\cdot 10^{-3}$ & $5.9\cdot 10^{-4}$ & $5.0\cdot 10^{-7}$ & $1.3\cdot 10^{-7}$\\
$128$  & $2.1\cdot 10^{-4}$ &  $2.4\cdot 10^{-4}$ & $3.8\cdot 10^{-3}$ & $1.4\cdot 10^{-4}$ & $1.2\cdot 10^{-7}$ & $5.9\cdot 10^{-8}$\\
\end{tabular}
}
\caption{Convergence results for the lubricated immersed boundary method applied to the problem of eccentric rotating cylinders.}
\label{tab:test2_lube}
\end{table}

\section{Applications}

\subsection{Wall-induced migration}
\begin{figure}[!ht]
         \begin{center}
         	\subfloat[]{\label{fig_onewall:exp}
	 \includegraphics[width=2.6 in]{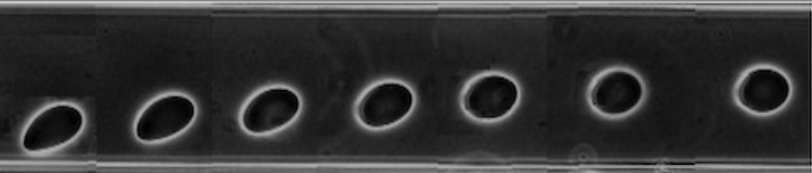}}\\
	 	\subfloat[]{\label{fig_onewall:sim_lube}
	 \includegraphics[width=3.2 in]{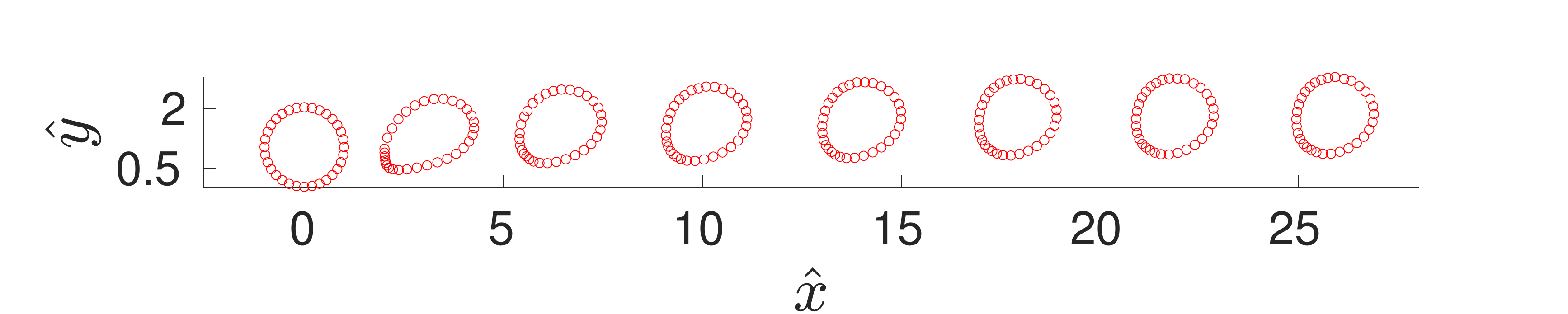}}\\
	 	\subfloat[]{\label{fig_onewall:sim_ib}
	 \includegraphics[width=3.2 in]{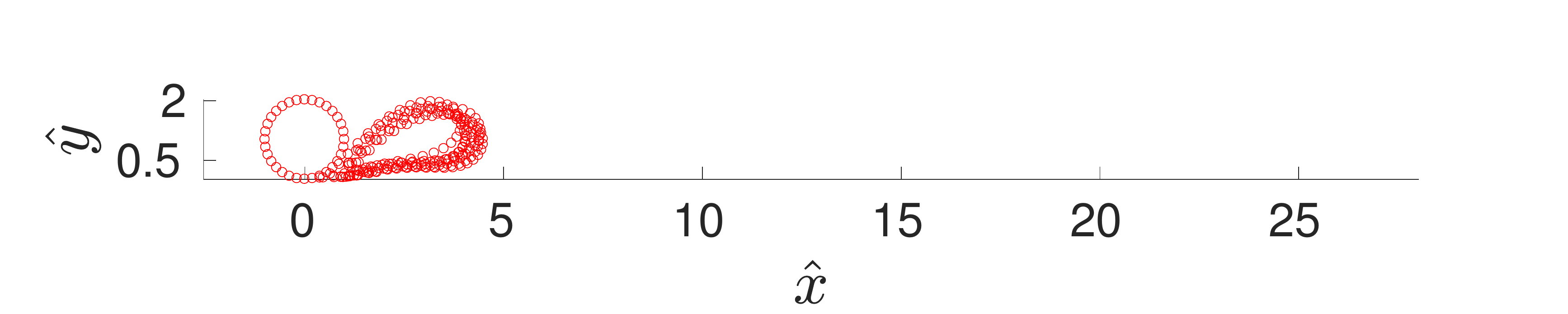}}
         \caption{\protect\subref{fig_onewall:exp} Experimental demonstration of wall-induced migration of an elastic vesicle, reprinted from \cite{coupier2008noninertial} with the permission of AIP Publishing, \protect\subref{fig_onewall:sim_lube} Wall-induced migration simulated by the lubricated immersed boundary method on a $64 \times 64$ grid. The position of the elastic ring is shown at several times as it moves through the channel and away from the wall. \protect\subref{fig_onewall:sim_ib} Analogous simulation by the standard immersed boundary method; the ring fails to migrate from the wall because of the underresolved relative velocity between the two boundaries. Axes are given in the dimensionless units $\hat{x} = x/R_0$, $\hat{y} = y/R_0$.}
         \label{fig_onewall}
         \end{center}
\end{figure}

Many fluid-structure interaction problems in biology, including blood flow, intracellular trafficking, and leaf aerodynamics involve the near-contact of multiple objects. Blood is a non-Newtonian fluid that exhibits several counter-intuitive behaviors, and one interesting such phenomenon is the emergence of a cell-free layer at the vessel walls in which no red blood cells are present. The cell-free layer has its origins in wall-induced migration, which causes deformable objects such as liquid droplets and red blood cells to acquire lift and move away from walls (see Figure \ref{fig_onewall}\protect\subref{fig_onewall:exp} and \cite{chaffey1965particle,uijttewaal1993droplet,cantat1999lift,coupier2008noninertial}). Because blood in the microcirculation is at low Reynolds numbers ($\text{Re} \approx 10^{-4}-10^{-2}$) so that the Stokes flow approximation is valid, directed motion away from the wall may seem to be at odds with the well-known scallop theorem that asserts no net displacement may be achieved by a time-symmetric motion in Stokes flow  \cite{purcell1977life}. The paradox is resolved in this case since symmetry is broken as the object deforms with the flow \cite{uijttewaal1993droplet,cantat1999lift}.

Simulating wall-induced migration by the standard immersed boundary method is challenging because of the issues discussed thusfar in resolving the lubricating layer that separates the deformable object from the wall. However, by applying the lubricated immersed boundary method, motion from the wall can be seen even on coarse grids. We have studied this phenomenon by simulating an elastic ring moving through a channel. The elastic ring is made up of immersed boundary points $\mb{X}_i$ for $i=0,\,1,\,\dots,\,N_q-1$ arranged in a loop and connected by springs having total energy
\begin{equation}
 E_\text{spring} = \frac{k_\text{spring}}{2} \sum_{i=0}^{N_q-1} \left(\frac{\norm{\mb{X}_{i+1}-\mb{X}_{i}}}{\Delta q}\right)^2\Delta q,
\end{equation}
and by torsional springs with total bending energy
\begin{equation}
E_\text{bend} = \frac{k_\text{bend}}{2} \sum_{i=0}^{N_q-1} \left(\frac{\norm{\mb{X}_{i-1}-2\mb{X}_{i}+\mb{X}_{i+1}}}{\Delta q^2} \right)^2\Delta q,
\end{equation}
where the indices are computed modulo $N_q$. We simulate the motion away from the wall using the standard and lubricated immersed boundary methods and compare the qualitative results (see Figure \ref{fig_onewall}). The elastic ring has radius $R_0 = \pi/4$, $k_\text{spring} = 1$, and $k_\text{bend} = 0.004$, and it is immersed in a periodic domain of unit cell size $L_x = 4 \pi$ by $L_y= 2 \pi$ that is filled with fluid having density $\rho = 1$ and viscosity $\mu = 0.2$. The viscosity is taken to be spatially uniform. This is the appropriate assumption for systems such as lipid-bilayer endosomes and red cell ghosts in which the internal and external fluids share the same material properties.

The channel walls are built by putting down a single line of tether points at $y = 0$. One wall of tether points is sufficient to define the channel because the periodic boundary conditions in $y$ imply that a particle will reach the wall regardless of whether it moves up or down.

Flow is established in the channel by applying a body force of the form $(f_x = \mu \sin(y),f_y=0)$. In the absence of an elastic ring, this establishes a unidirectional flow $(u = \sin(y),v=0)$ in the domain. We use a $64 \times 64$ grid (in which case $\Delta x \neq \Delta y$) and timestep $\Delta t = 0.005 \Delta x$, and 32 discrete points around the circumference of the ring. The elastic ring is initialized nearby the wall with its center of mass at $(x,y) = (0, R_0+\pi/192)$. This implies that at the beginning of our simulation the ring is at a distance of $\Delta y/6$ to the wall. The simulations are run up to a time of $T = 28$.

In these simulations of elastic vesicles, we use a modified fluid solver that has been used previously to simulate inkjet printers \cite{yu2003coupled} and is based on the finite element projection of \cite{almgren1996numerical}. A marker-and-cell (MAC) fluid velocity is constructed at each time step in which the horizontal and vertical velocities are defined at cell edges and are not co-located. \edone{The advantage of this fluid solver is that the MAC velocity is useful as an interpolating field since advected elastic structures tend to suffer less leakage and conserve volume better than those advected on collocated grids \cite{griffith2012volume}.}

Figure \ref{fig_onewall} shows that, at a qualitative level, the lubricated immersed boundary method is able to capture the sequence of deformations during wall-induced migration. The elastic ring simulated by the lubricated immersed boundary method deforms and moves away from the wall, whereas the ring simulated by the standard immersed boundary method becomes stuck and cannot travel with the flow (see Movies \texttt{onewall\_ib64.gif} and \texttt{onewall\_lube64.gif}).

\subsection{Channel flow}
\begin{figure}[!htp]
\captionsetup[subfigure]{labelformat=empty}
\vspace{-1.8in}

         	 \hspace{-1.2in} \subfloat[]{\label{fig_twowall:ib1}
	\includegraphics[width=3.2 in]{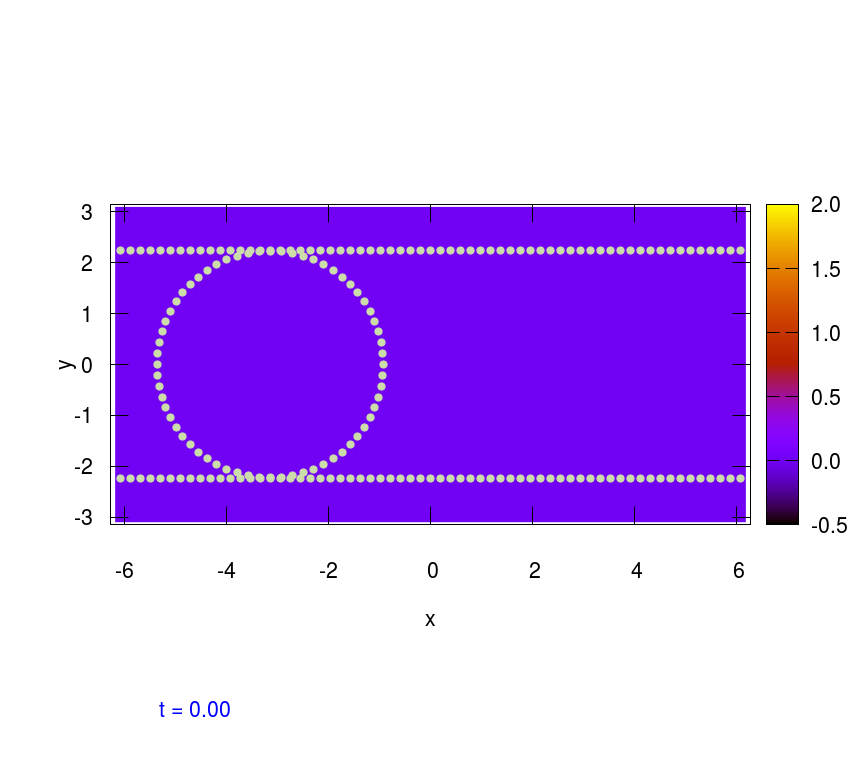}}
	       \subfloat[]{\label{fig_twowall:lube1}
	 \includegraphics[width=3.2 in]{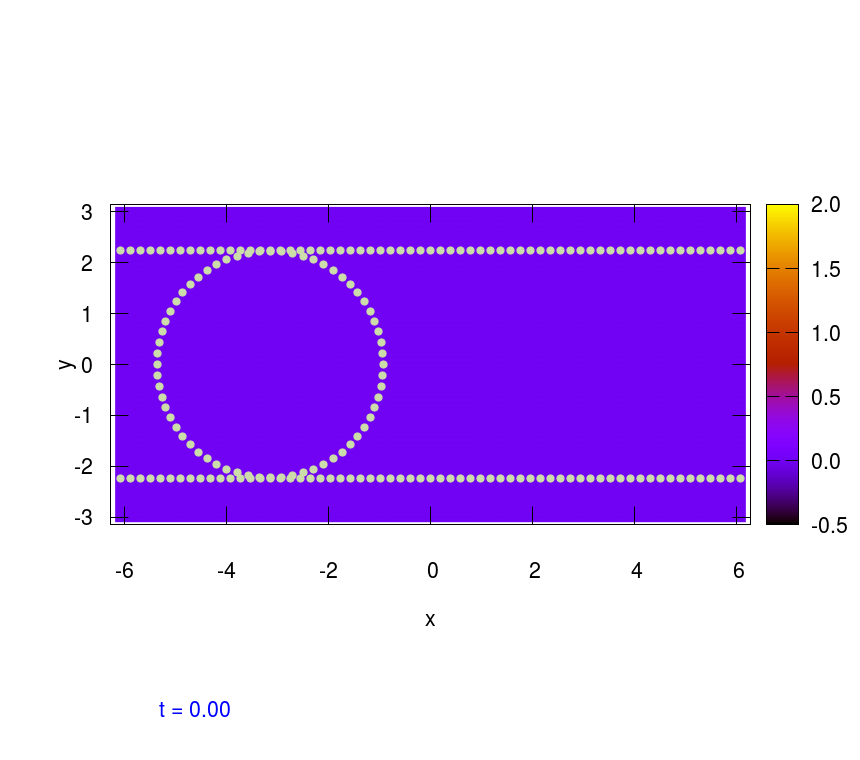}}\\
	 \vspace{-0.7in}
	 
         	\hspace{-1.2in} \subfloat[]{\label{fig_twowall:ib2}
	 \includegraphics[width=3.2 in]{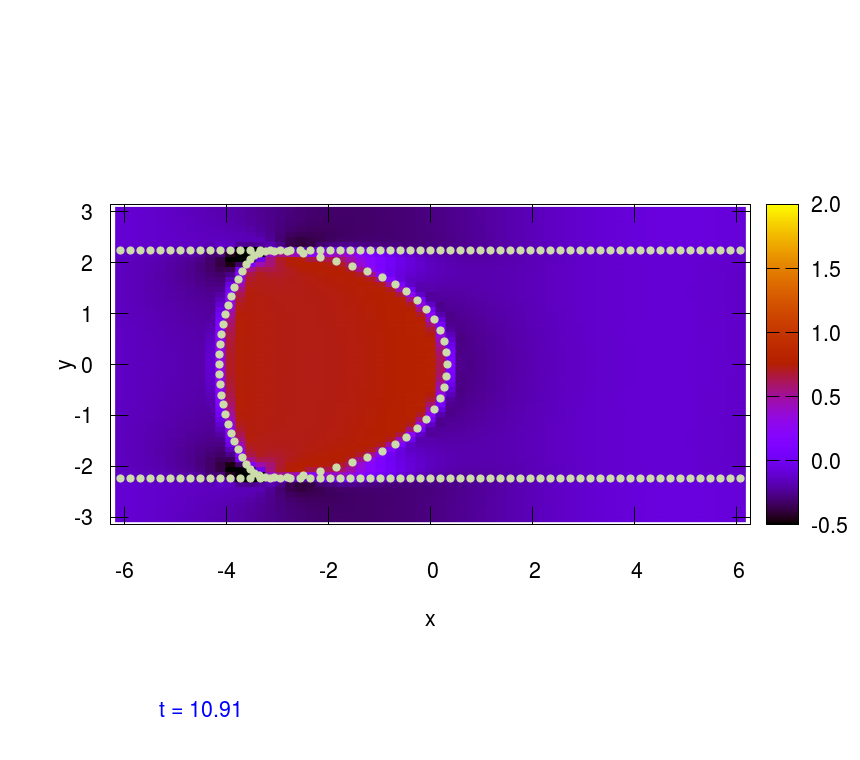}}
	 	\subfloat[]{\label{fig_twowall:lube2}
	 \includegraphics[width=3.2 in]{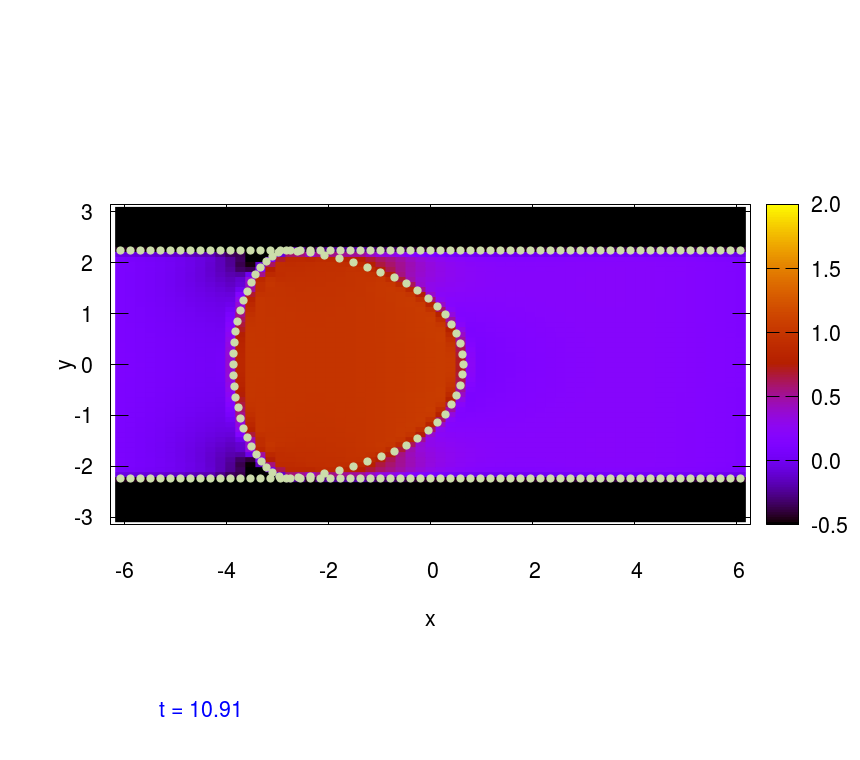}}\\
	 	 \vspace{-0.7in}
	 	 
	        \hspace{-1.2in} \subfloat[]{\label{fig_twowall:ib4}
	 \includegraphics[width=3.2 in]{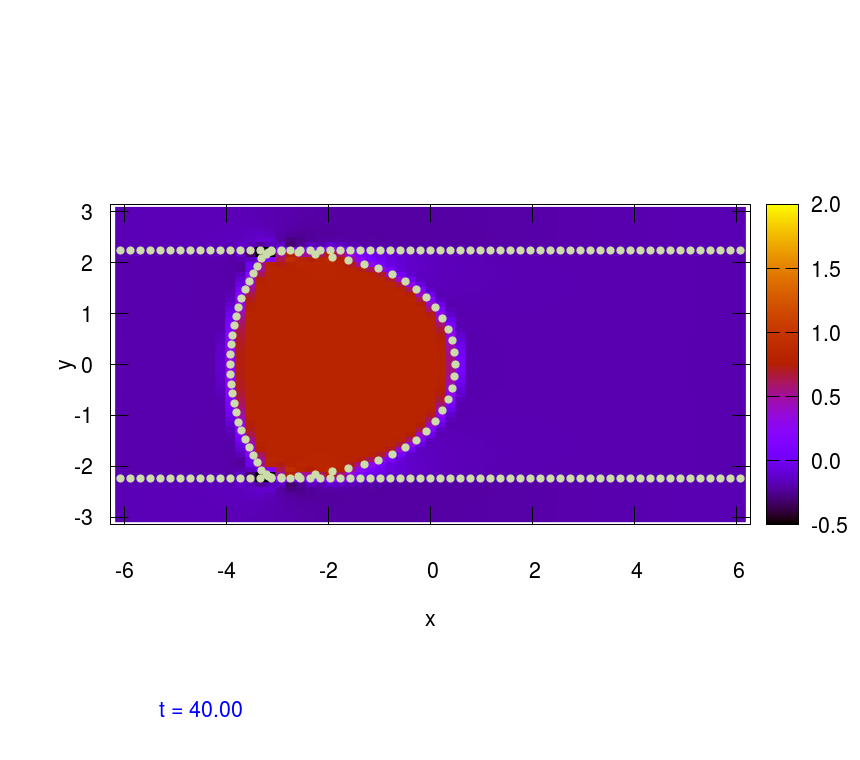}}
	 	\subfloat[]{\label{fig_twowall:lube4}
	 \includegraphics[width=3.2 in]{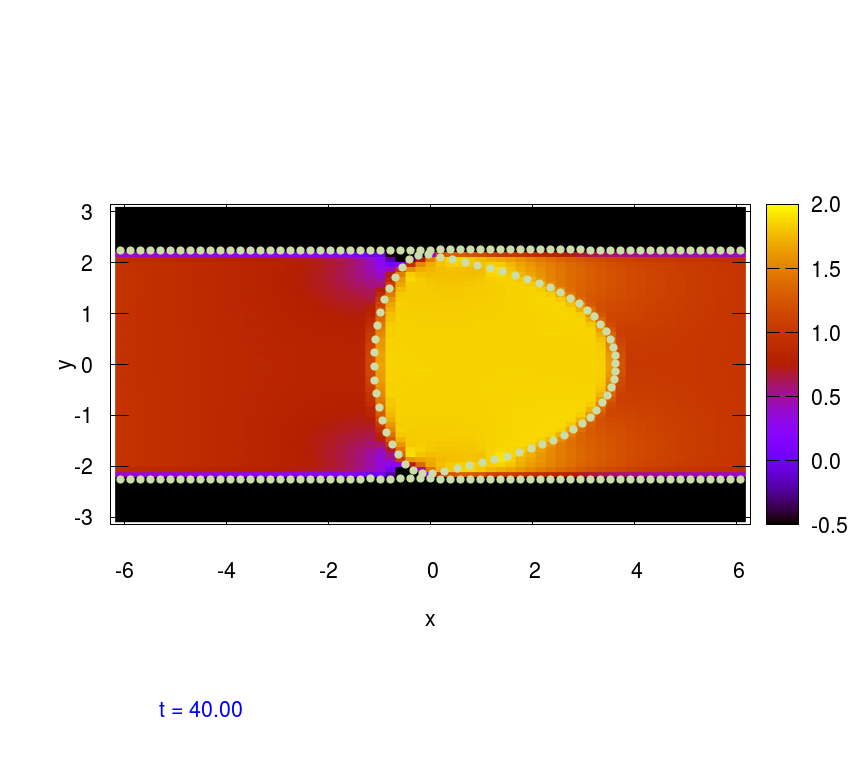}}\\
         \caption{Comparison between simulations of channel flow using the standard (left column) and lubricated (right column) immersed boundary methods, with the color map representing the pressure. An elastic ring is initialized at a distance of $\Delta y/6$ from both walls and a uniform force is applied to the ring. In the standard immersed boundary method simulation, the ring sticks to the wall near its starting position, whereas in the lubricated immersed boundary method simulation the ring separates from the channel walls and moves in the direction of the applied force.}
         \label{fig_channel}
\end{figure}

\begin{figure}[!htp]
         \begin{center}
	         \subfloat[$u$, standard method]{
	 \includegraphics[width=5 in]{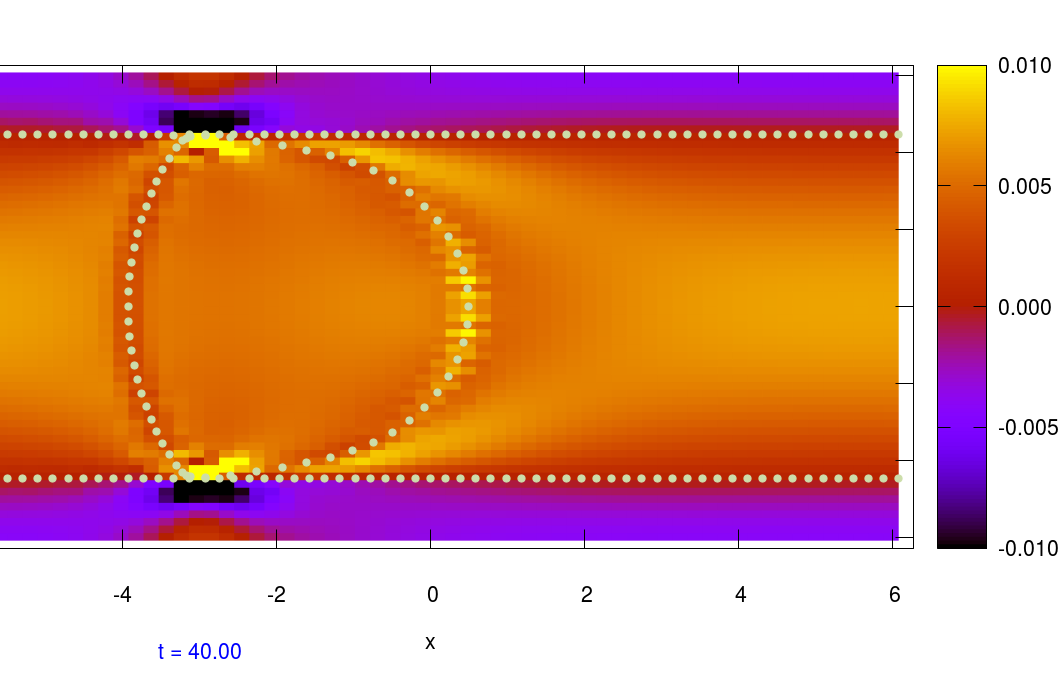}}\\
	 	 	 \vspace{-0.1in}
	 	\subfloat[$u$, lubricated method]{
	 \includegraphics[width=5 in]{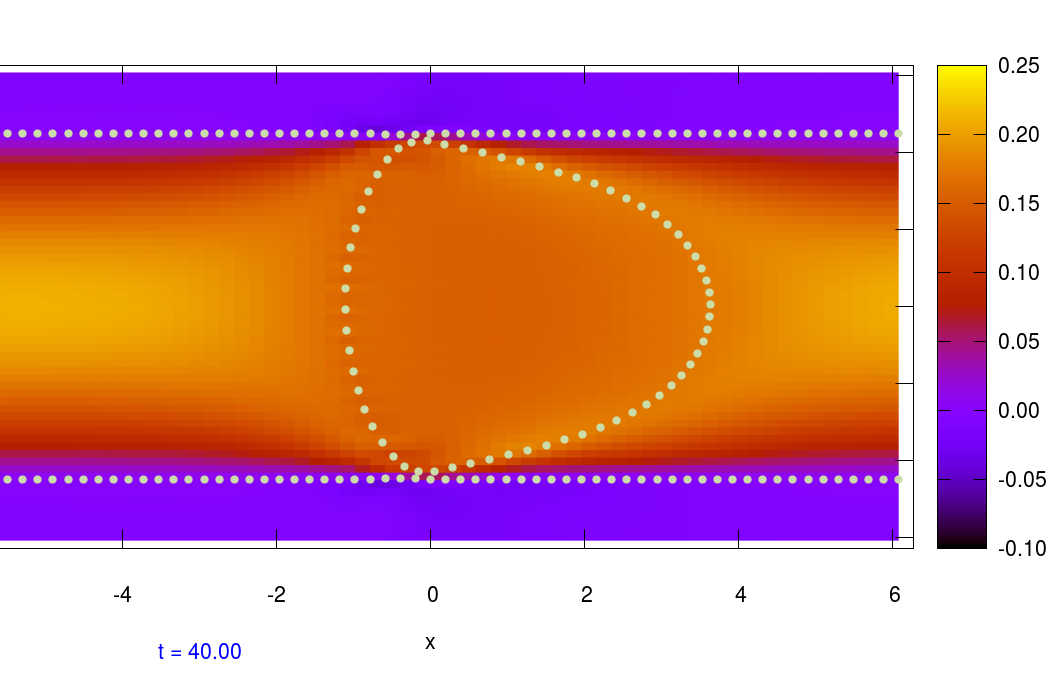}}\\
         \caption{Comparison between simulations of channel flow using the standard and lubricated immersed boundary methods, with the color map representing horizontal velocities at time $t = 40$. There is an order of magnitude difference in the horizontal velocities obtained; in the standard immersed boundary method simulation, the channel is effectively clogged.}
         \label{fig_channel_u}
         \end{center}
\end{figure}

Finally, we perform a simulation of an elastic ring flowing through a narrow channel. This example is inspired by the fluid-structure interaction that takes place in the microcirculation when red blood cells squeeze through capillaries as small as a few \textmu m in diameter.

The size of the elastic ring in the following simulation is chosen so that it fits just within a channel defined by two walls of tether points. We use an elastic ring of radius $R_0 = 9\pi^2/40$ that is initialized at a distance of $\pi/192$ from each wall. As in the previous section, we use a $64 \times 64$ grid on a domain of size $L_x = 4 \pi$ by $L_y= 2 \pi$ so that, in terms of the grid spacing, the ring is initialized at a minimal distance of $\Delta y/6$ from each wall. The ring is forced in the $+x$-direction by applying a uniform force to each immersed boundary point. Here we discretize the ring with $64$ immersed boundary points around its perimeter and use the timestep $\Delta t = 0.01 \Delta x$. The simulation is run up to $T = 40$. Note that for this problem, we must compute two height functions, one for the distance of the vesicle to each wall. 

When the lubricated method is used, the ring flows along in the channel, adopting bullet-like shapes reminiscent of red blood cells in capillaries. This is in contrast to results obtained by the standard immersed boundary method, in which the ring becomes stuck at its starting point (Figure \ref{fig_channel}). The horizontal velocity field computed by the standard method displays grid artifacts in the velocities in the narrow gap separating the elastic ring from the walls, and small velocities throughout the rest of the channel, indicating a clog. The horizontal velocities computed by the lubricated immersed boundary method, on the other hand, do not display these grid artifacts in the lubrication layers, and are over an order of magnitude greater in the rest of the channel (Figure \ref{fig_channel_u}). See Movies \texttt{channel\_ib64\_p.gif} and \texttt{channel\_lube64\_p.gif}, which show the evolution of the pressure field, and Movies \texttt{channel\_ib64\_u.gif} and \texttt{channel\_lube64\_u.gif}, which show the evolution of the horizontal velocity, for a comparison of the results obtained by the standard and lubricated methods.

\section{Conclusions}
We have proposed an immersed boundary method that uses aspects of lubrication theory to resolve thin fluid layers between elastic structures. The two key ingredients of the method are (i) the observation that across a lubrication layer the normal derivatives of velocity at the boundaries determine the relative velocity, and (ii) the well-established jump conditions that relate these normal derivative of velocity to the Lagrangian force on the boundary. We have applied the lubricated immersed boundary method to problems of increasing complexity and have found that the method outperforms the standard immersed boundary method when the lubrication layer is small compared to the fluid grid spacing. The lubricated method makes it possible on coarse fluid grids to capture interesting phenomena, such as the lift and reverse flow in the problem of eccentric rotating cylinders, and wall-induced migration in the problem of an elastic vesicle near a wall in shear flow.

Note that, although in the problems considered here the immersed boundaries were initialized to be close to touching, immersed boundaries that are initially well-separated can also come within subgrid distances over the course of a simulation. This effect is particularly apparent in underresolved 3D computations with significant differences in velocities from gridpoint to gridpoint, as occurs for instance in simulations of red cells moving through capillaries \cite{vc_ib} and flagellar bundling in E. Coli \cite{lim2012}. As noted in the Introduction, in this last case initially separated immersed boundaries have even been observed to cross, making it a good candidate to which to apply the lubricated method.

Although we have demonstrated the effectiveness of the present formulation in several settings, much work remains to be done to make the lubricated immersed boundary method readily applicable to real-world problems. The results presented here are restricted to 2D, and although much of the mathematical formulation extends to 3D, it is not clear to what extent the specific numerical approach (e.g.~the definition of the height function in terms of piecewise cubic curves) will extend to the challenges of 3D problems. Further investigation may indicate that alternative formulations are preferable to the one presented here; for example, the current formulation is not symmetric with respect to upper and lower surfaces, a property that would be desirable for some applications. \edone{(However, for those applications in which one surface tends to be significantly flatter than the other, an asymmetric formulation may be preferable.)} \edtwo{Moreover, in many applications it may be more convenient to impose physical boundary conditions within the fluid solver rather than using a tethered surface.} In future work, we plan to find appropriate lubrication corrections in these types of situations as well.

\section{Acknowledgments}
  
We thank Charles Peskin for useful conversations \edthree{and the anonymous reviewers for their helpful feedback}. TGF acknowledges funding through National Science Foundation grant DMS-1502851. CHR was supported by the Applied Mathematics Program of the U.S. Department of Energy (DOE) Office of Advanced Scientific Computing Research under contract DE-AC02-05CH11231.
  
  \appendix
  \section{Computing the height on the upper surface}
  \label{app:height}
  
\begin{figure}[!ht]
         \begin{center}
	 \includegraphics[width=3 in]{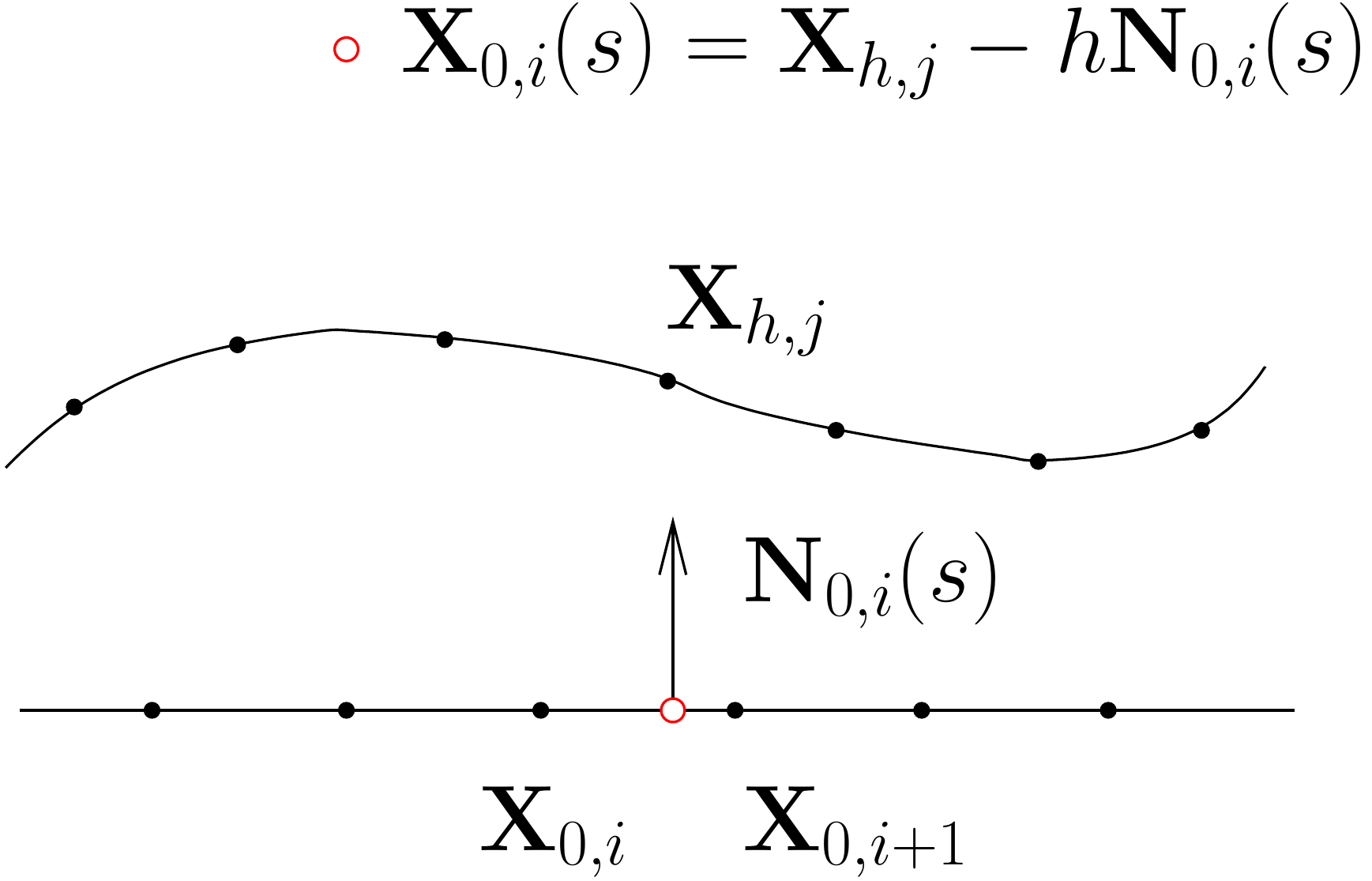}
         \caption{\edone{Computing the height $h$ for a point $\mb{X}_{h,j}$ on the discrete upper surface. The corresponding point on the lower surface is denoted by $\mb{X}_{0,i}(s)$, where $\mb{X}_{0,i}(s) = \mb{a}_0+\mb{a}_1 s+\mb{a}_2 s^2 +\mb{a}_3 s^3$ for $s \in [0,1]$. That is, the point $\mb{X}_{0,i}(s)$ is located on a section of the lower surface defined by a cubic curve between the discrete points $\mb{X}_{0,i}$ and $\mb{X}_{0,i+1}$ The height is defined with respect to the lower surface normal $\mb{N}_{0,i}(s)$ at the point $\mb{X}_{0,i}(s)$.}}
         \label{fig:height_upper}
         \end{center}
\end{figure}

Consider a point $\mb{X}_{h,j}$ on the upper surface (Figure \ref{fig:height_upper}). We approximate the lower surface by a piecewise cubic curve as before. That is, given a section of the lower surface with endpoints $\mb{X}_{0,i}$ and $\mb{X}_{0,i+1}$, the curve between these endpoints is $\mb{X}_{0,i}(s) = \mb{a}_0+\mb{a}_1 s+\mb{a}_2 s^2 +\mb{a}_3 s^3$ for $s \in [0,1]$, where the coefficients are defined in an analogous manner to \eqref{eq:cub1}--\eqref{eq:cub4}.
We wish to find the pair $(s,h)$ and the index $i$ such that
\begin{equation}
  \mb{X}_{0,i}(s)+h \mb{N}_{0,i}(s) = \mb{X}_{h,j}, \label{eq:upper_height}
\end{equation}
where as before $\mb{N}_{0,i}(s)$ is defined by $\mb{N} =  (-T^2,T^1)/\norm{\mb{T}}$ through the tangent vector, which in this case is given by differentiating the cubic surface, i.e.~
\begin{equation}
\mb{T}_{0,i}(s) = \frac{\ud \mb{X}_{0,i}(s)}{\ud s} = \mb{a}_1+2 \mb{a}_2 s +3 \mb{a}_3 s^2.
\end{equation}
Since the normal vector itself depends on the unknown $s$, the formulation \eqref{eq:upper_height} is subtly different from the analogous formulation \eqref{eq:lower_height} of the height above a point on the lower surface. To solve for $s$, we take the dot product of \eqref{eq:upper_height} with $\mb{T}_{0,i}(s)$ to obtain the 5\textsuperscript{th} order polynomial equation
\begin{align}
  &\,\mb{X}_{0,i}(s)\cdot \mb{T}_{0,i}(s) = \mb{X}_{h,i}\cdot \mb{T}_{h,i}(s) \notag \\
  \iff &\left(\mb{a}_0+\mb{a}_1 s+\mb{a}_2 s^2 +\mb{a}_3 s^3-\mb{X}_{h,j}\right)\cdot \left(\mb{a}_1+2 \mb{a}_2 s +3 \mb{a}_3 s^2\right) = 0
\end{align}
and once more use Newton's method to compute the solution. Given $s$, it is straightforward to solve for the height via $h_j = \left(\mb{X}_{h,j}-\mb{X}_{0,i}(s)\right) \cdot \mb{N}_{0,i}(s)$. As before, it is important to have a good initial guess to ensure the efficient convergence of Newton's method. To obtain this initial guess, we approximate the lower surface by a piecewise linear function so that $\mb{X}_{0,i}(s) = \mb{X}_{0,i}+s\left(\mb{X}_{0,i+1}-\mb{X}_{0,i}\right)$, but we use the continuously varying normal vector
\begin{equation}
\mb{N}_{0,i}(s) = \mb{N}_{0,i}+s\left(\mb{N}_{0,i+1}-\mb{N}_{0,i}\right).
\end{equation}
The continuously-varying normal is required to ensure there is a region of finite width around the lower surface above which heights can be computed. Now, we wish to find $s_0$ such that
\begin{align}
 \mb{X}_{h,j}&=\mb{X}_{0,i}(s_0)+h \mb{N}_{0,i}(s_0) \notag \\ &=\mb{X}_{0,i}+s_0\left(\mb{X}_{0,i+1}-\mb{X}_{0,i}\right)+h\left(\mb{N}_{0,i}+s_0\left(\mb{N}_{0,i+1}-\mb{N}_{0,i}\right)\right).
\label{eq:upper_height_guess}
\end{align}
Taking dot products with $\mb{T}_{0,i}(s_0)$ yields the quadratic equation
\begin{equation}
a s^2+b s+c=0,
\end{equation}
where
\begin{align}
 a &=\left(\mb{X}_{0,i+1}-\mb{X}_{0,i}\right)\cdot\left(\mb{T}_{0,i+1}-\mb{T}_{0,i}\right)\\
 b &=\left(\mb{X}_{0,i+1}-\mb{X}_{0,i}\right)\cdot\mb{T}_{0,i}+\left(\mb{X}_{0,i}-\mb{X}_{h,j}\right)\cdot\left(\mb{T}_{0,i+1}-\mb{T}_{0,i}\right)\\
 c &= \left(\mb{X}_{0,i}-\mb{X}_{h,j}\right)\cdot\mb{T}_{0,i}.
\end{align}
Of the two roots $(-b\pm\sqrt{b^2-4ac})/(2a)$, the negative root is the relevant one since it remains finite in the limit of zero curvature in which $a \to 0$. We have found that using the quadratic formula to actually compute the roots can lead to a loss of precision; instead we use the numerically robust method described in \cite{press2007numerical}. 


\bibliography{mybiblio}

\end{document}